\documentclass{amsart}
\usepackage{graphicx} 
\usepackage{enumerate}
\usepackage{enumitem}

\usepackage{xcolor}
\usepackage{graphicx}
\usepackage{epsfig}
\usepackage{tikz}

\usepackage{latexsym}
\usepackage{amssymb}
\usepackage{amsmath}
\usepackage{amsthm}
\usepackage{amscd}
\usepackage{url}

\usepackage{enumerate}

\theoremstyle{plain}
\newtheorem{thm}{Theorem}[section]
\newtheorem{lem}[thm]{Lemma}
\newtheorem{defin}[thm]{Definition}
\newtheorem{prop}[thm]{Proposition}

\theoremstyle{definition}
\newtheorem{rem}{Remark}

\newtheorem{ex}{Example}

\newtheorem{cor}[thm]{Corollary}

\newcommand{\re}{\operatorname{Re}}
\newcommand{\im}{\operatorname{Im}}
\newcommand{\Bal}{{\mathrm{Bal}}}
\newcommand{\Taut}{{\mathrm{Taut}}}
\newcommand{\Gal}{{\mathrm{Gal}}}

\newcommand{\veigenvalue}{{virtual eigenvalue}}
\newcommand{\veigenvalues}{{virtual eigenvalues}}

\newcommand{\uzero}{\epsilon}
\newcommand{\um}{\frac{10 \uzero}{\kappa}}

\newcommand{\R}{\mathbb{R}}
\newcommand{\C}{\mathbb{C}}
\newcommand{\Z}{\mathbb{Z}}
\newcommand{\N}{\mathbb{N}}
\newcommand{\Q}{\mathbb{Q}}
\newcommand{\T}{\mathbb{T}}

\newcommand{\Hol}{\mathrm{Hol}}

\author{Francisco Arana--Herrera}
\address{Department of Mathematics, University of Maryland, College Park, MD USA}
\email{farana@umd.edu}

\author{Jon Chaika}
\address{Department of Mathematics, University of Utah, Salt Lake City, UT USA}
\email{chaika@math.utah.edu}

\author{Giovanni Forni}
\address{Department of Mathematics, University of Maryland, College Park, MD USA and Laboratoire AGM, CY Cergy Paris Universit\'e, France}
\email{gforni@umd.edu}

\title{Weak mixing in rational billiards}

\begin{document}

\begin{abstract}
    We completely characterize rational polygons whose billiard flow is weakly mixing in almost every direction as those which are not almost integrable, in the terminology of Gutkin~\cite{Gu86}, modulo some low complexity exceptions. This proves a longstanding conjecture of Gutkin~\cite{Gu84}. This result is derived from a complete characterization of translation surfaces that are weakly mixing in almost every direction: they are those that do not admit an affine factor map to the circle. 
\end{abstract}

\maketitle

\section{Introduction}

\subsection{Background and motivation} Planar billiards are an important and well-studied class of dynamical systems, in particular, because they display, at least in a simplified form, the main features of the different classes of dynamical behaviours one can witness in nature; see for instance the survey \cite{Kat_05}. Billiards in polygons are a prime example of {\it parabolic dynamical systems}, characterized by slow, i.e., polynomial, divergence of nearby orbits with respect to time. These features have also attracted the interest of many mathematical physicists studying different quantum aspects of dynamical systems; see for instance \cite{RB81}.

The dynamics of billiards in polygons presents a fundamental dichotomy between the rational case, in which the phase space is foliated by invariant surfaces 
and the non-rational case \cite{MT02}. Rational billiards belong to a class of systems which have been called {\it pseudo-integrable} \cite{RB81}, in analogy with the notion of
an integrable system in classical mechanics, where most of the phase space is foliated by invariant tori. 

This paper is concerned with the ergodic theory of billiards in rational polygons. The main result is a complete characterization of which rational polygons have a billiard flow which is weakly mixing in almost every direction. We recall that  A.~Katok \cite{Kat80} proved that Interval Exchange Transformations (IETs)  and linear flows on translation surfaces are never mixing, hence it is natural to see what billiards are weakly mixing. 
Incidentally, together with A. Katok's absence of mixing result, our weak mixing theorem confirms that, like the typical measure preserving transformation, the typical billiard flow is weakly mixing but not strongly mixing.

The study of rational billiards, and of geodesic flows on rational polyhedra, has a long history. The topological properties of billiard flows on rational polyhedra were first investigated in \cite{FK36} and later in \cite{KZ75,BKM}. The ergodic theory of rational billiards on Lebesgue almost every invariant surface had to wait until the introduction by G.~Rauzy \cite{Ra79}, W.~Veech \cite{Ve82}, and H.~Masur \cite{Ma82} of a renormalization scheme for IETs and their suspensions with respect to piecewise constant roof functions, commonly known as translation surfaces. Renormalization has proved to be an extremely powerful tool in the study of ergodic properties for Lebesgue almost every element of the classes of systems mentioned above. In particular, it has lead to the 
proof of unique ergodicity for Lebesgue almost every IET \cite{Ve82,Ma82}. However, these results do not apply
to rational polygons, since the latter correspond to an exceptional, i.e., zero Lebesgue measure, set in the space of all IETs or translation surfaces. Unique ergodicity for Lebesgue almost every invariant surface of a given rational polygon was nevertheless proved
a few years later in \cite{KMS} using the renormalization dynamics of the Teichm\"uller geodesic flow, introduced in the work of H.~Masur \cite{Ma82}, and the non-divergence property of the so-called Teichm\"uller horocycle flow, introduced by H.~Masur in \cite{Ma85}. 

Results on the conjectural weak mixing properties of IETs
and translation surfaces, which will be recalled below, have 
so far not been extended to billiards in rational polygons, with the exception of special cases including regular polygons with at least $5$ edges \cite{AvDe}. Contrary to unique ergodicity, there exist counterexamples to weak mixing in almost every direction in the wider class of translation surfaces; for instance, any branched cover of the flat torus provides such a counterexample. Our results provide, in particular and for the first time, a complete classification of which translation surfaces are weak mixing in almost every direction. The crucial breakthroughs on which our work is based are the measure rigidity and orbit closure results of A.~Eskin and M. Mirzakhani \cite{EM}, and A.~Eskin, M.~Mirzakhani and A.~Mohammadi \cite{EMM}, together with
important follow-ups in work of the second author and 
A.~Eskin \cite{CE15}, as well as the work of 
S.~Filip \cite{mono} on the dynamics of the Kontsevich--Zorich cocycle \cite{Ko97}, \cite{KZ97}, \cite{F02}, \cite{F06}, \cite{F11}.

\subsection{Statement of results}
A translation surface is a pair $(X,\omega)$, where $X$ is a Riemann surface and $\omega$ is a non-zero holomorphic 1-form. The holomorphic 1-form $\omega$ induces a singular translation structure on $X$, giving, in particular, a locally Euclidean metric on $X \setminus \Sigma$, where $\Sigma$ is the finite set of zeros $\omega$. These zeroes correspond to cone type singularities of the metric with cone angles in $2 \pi \N$. From this structure we obtain an area form on $X$ and a flow $F_\theta^t$ in each direction $\theta \in [0,2\pi)$, that preserves the area. Note that, because of the singularities, each flow is defined only on a set of full measure with respect to the area form. As mentioned above, these dynamical systems, i.e., flows on translation surfaces, are much studied and have explicit connections to billiards in polygons, smooth flows on surfaces, and IETs. For every translation surface, the flow in almost every direction is uniquely ergodic \cite{KMS}. While some translation surfaces are not weakly mixing in any direction, for the typical translation surface, the flow in the vertical direction is weakly mixing \cite{AvFo}. In this paper we completely classify which translation surfaces have the property that the flow in almost every direction is weakly mixing. More concretely:

\begin{thm} 
\label{thm:main} Let $(X,\omega)$ be a translation surface and $F_\theta$ denote its flow in the direction $\theta \in [0,2\pi)$. Then, the following are equivalent:
\begin{enumerate}
\item The flow $F_\theta$ is weakly mixing for almost every $\theta \in [0,2\pi)$.
\item The flow $F_\theta$ is weakly mixing for some $\theta \in [0,2\pi)$.
\item There exists no non-zero integer absolute cohomology class in the tautological plane of $(X,\omega)$, i.e., the plane in $H^1(X;\mathbb{R})$ spanned by classes of the real and imaginary parts of $\omega$. 
\end{enumerate}
\end{thm}

From the above result we derive a complete characterization of rational polygons whose billiard flow is weakly mixing in almost every invariant surface, i.e., in almost every direction, thereby confirming a strong version of a conjecture of 
E.~Gutkin \cite{Gu84}. 
Indeed, Gutkin introduced a class of polyhedral surfaces and stated that  ``It is reasonable to expect that for almost integrable polyhedral surfaces outside
of this class the billiard flow is typically weakly mixing''; see \cite{Gu84}, page 570. Note that the terminology of
\cite{Gu84} was abandoned already in \cite{Gu86} and, while
rational billiards correspond via doubling to ``almost integrable polyhedral surfaces'', the term ``almost integrable'' was restricted to the special class he had introduced.  

An {\it almost integrable} polygon is a polygon drawn on the ``grid'', or ``lattice'' in the terminology of Gutkin, spanned by reflecting one of the completely integrable polygons along its sides ad infinitum.
The completely integrable polygons are all rectangles and
the $(\pi/2, \pi/4, \pi/4)$, $(\pi/3, \pi/3, \pi/3)$, and $(\pi/2, \pi/3, \pi/6)$ triangles; see \cite[\S 1.4]{MT02}.

\smallskip
In this paper, as a consequence of Theorem \ref{thm:main}, we prove the following complete classification of the polygons whose billiard flow is weak mixing in almost every direction.

\begin{cor} 
\label{cor:main}
The billiard flow in a rational polygon $P$ is weakly mixing on almost every invariant surfaces, i.e., in almost every
direction, unless $P$ is an almost integrable polygon or all of its angles are in $\{\frac \pi 2 ,\frac {3\pi}2\} $ and if oriented to have horizontal/vertical sides, all the horizontal lengths or all the vertical lengths are commensurable.\footnote{The horizontal/vertical side lengths are commensurable if all their ratios are rational.}  
\end{cor}

\subsection{Related results}
Weak mixing for almost every dynamical  system in a given class was proved for certain classes of IETs earlier than for flows on translation surfaces. Veech \cite{Vmetric} proved that weak mixing was typical for IETs with \emph{type W permutations}. 
M.~Boshernitzan and A.~Nogueira \cite{BoNo} gave a different proof that  provided an explicit, natural, full measure set. Nogueira and D.~Rudolph \cite{NoRu} showed that the typical non-rotation-class IET is topologically weakly mixing. A.~Avila and the third author \cite{AvFo} then proved that the typical non-rotation-type IET and the vertical flow on a typical translation surface of genus at least $2$ are weakly mixing. Avila and M.~Leguil \cite{AvLe} strengthened this to show that in these situations, weak mixing held off of sets of positive Hausdorff codimension. In the work closest to this paper, Avila and V.~Delecroix \cite{AvDe} classified which Veech surfaces are weak mixing in almost every direction, hence in particular they established weak mixing in the typical direction for the billiard flow in a regular polygon with at least $5$ edges. More recently \emph{effective weak mixing} has been a topic of extensive study \cite{BS14, BS18, BS20, BS21, BS22, Fo22, AFS23}.

The ergodic theory of billiards in non-rational polygons is
much less understood. 
S.~Kerckhoff, H.~Masur and J.~Smillie proved \cite{KMS} that there exists a $G_\delta$-dense set of polygons whose flow is ergodic on the full unit tangent bundle. This result is based on a fast approximation argument which leverages the unique ergodicity of rational polygons on almost all invariant surfaces.  There is an effective version given later in \cite{Vo97} which requires super-exponential approximation of the polygon by rational ones. More recently, the last two authors \cite{CF} were able to prove that there exists a $G_\delta$-dense set of polygons whose flow is weakly mixing on the full unit tangent bundle by an approximation argument which leverages the equivalence of weak mixing with square ergodicity.

\subsection{Sketch of proof}
As is common in the study of weak mixing, we use the Veech criterion \cite{Vmetric} as the main technical tool in our proof. The basic idea of this criterion is that, if $u$ is a measurable eigenfunction for a flow $F$ on a space $X$, so that $u \circ F^s =e^{2 \pi i \alpha s} \cdot u$ for some $\alpha \in \mathbb{R}$, then, by Lusin's theorem, if $c>0$, $\epsilon>0$ is small enough, and there exists a time $r \in \mathbb{R}$ so that 
\begin{equation}\label{eq:wmix}\mathrm{Leb}(\{x \in X: d(F^rx,x)<\epsilon\})>c,
\end{equation} 
then we get the following restriction on the position of $\alpha$:
\[
d_\R(r \cdot \alpha,\Z) \ll 1
\]
Indeed, by Lusin's theorem, $u(F^r x)=e^{2\pi i r\alpha} \cdot u(x)$ needs to be close to $u(x)$ on a set of $x$'s of large measure where $u$ is continuous, so $e^{2 \pi i r\alpha}$ is close to $1$. Under reasonable compactness assumptions, for the vertical flow on a translation surface $(X,\omega)$ of genus $g\geq 1$, if $g_t$ denotes the Teichmüller geodesic flow and $X_t$ denotes the Riemann surface underlying $g_t \omega$, a vector ${\vec{ r}}=(r_1, \dots, r_{2g})$ whose entries satisfy \eqref{eq:wmix} is precisely  given by $g_t \cdot \im (\omega) \in H^1(X_t, \R) \approx \mathbb{R}^{2g}$, i.e., by the action of the \emph{Kontsevich-Zorich cocycle} \cite{Ve82}, \cite{Zo97}, \cite{Ko97}, \cite{KZ97} on the vector $\im (\omega)\in H^1(X, \R)$. The standard strategy to prove weak mixing is based on the fact that the Kontsevich--Zorich cocycle is isometric \cite{Vmetric} or  typically expands vectors \cite{F02}, \cite{AV07}; in this case the relevant vector is $\alpha \im (\omega)$. This allows one to apply the Veech criterion to rule out that $\alpha$ is an eigenvalue \cite{AvFo}. 

Complications arise because we only know that $d_\R(r \cdot \alpha,\Z)$ is small. We need to consider the difference $\alpha \im (\omega)-\vec{n}$ for an integer vector $\vec{n} \in H^1(X, \Z)\approx \Z^{2g}$ and allow $\alpha$ to vary in $\R\setminus \{0\}$. Thus, in order to apply this strategy, one seeks to show that, given a small interval of ``\veigenvalues'', i.e., a small interval $J \subset \R\setminus \{0\}$ and $\vec{n}\in H^1(X;\mathbb{Z})$ so that for all $\alpha \in J$ we have
\[
d\left(g_t \cdot \alpha \im (\omega),g_t \cdot \vec{n} \right)<\delta,
\]
then, for every such $\alpha$, there exists $s_{\alpha}>0$ so that 
\[
d\left(g_{t+s_\alpha} \cdot \alpha \im(\omega), H^1(X_t,\Z) \right)>\delta
\]
and $g_{t+s_\alpha}\omega$ is in a fixed compact set of moduli space.\footnote {Here $d$ can be taken to be any fiberwise metric on the cohomology bundle with fiber $H^1(X,\R)$ at every Riemann surface $X$, and it is most natural to consider $d = d_{X}$ to be the Hodge metric induced by the complex structure of $X$ on $H^1(X;\R)$.}

The naive hope would be that one can control the growth of $J$ as it evolves through the cocyle and guarantee there is a single time $s=s(J) > 0$ that works simultaneously for all $\alpha \in J$. Indeed, if $J_t := \{g_t \cdot \alpha \im(\omega):\alpha \in J\}$ and we assume 
\begin{itemize}
\item  $\vert g_{s} \cdot J_t \vert<\frac 1 {100}$,
\item $
d(g_{t+s} \cdot \beta \im (\omega),H^1(X_{t+s};\Z)) > 1/50$ for some $\beta \in J$,
\item and $g_{t+s}\omega $ is in a compact set, 
\end{itemize} 
then this is achieved. Unfortunately, this is often not possible, for instance, if the times at which the first two bullet points are simultaneously satisfied the third bullet point is never satisfied. This issue leads to our single small interval $J_t$ splitting into many subintervals. Controlling this splitting phenomena is a cause for much of the difficulty in our arguments. 

Many of our arguments rely on work of the last two named authors \cite{CF}, who proved strong expansion results in the symplectic complement of the tautological plane $\mathrm{Taut}(X,\omega):= \mathrm{span}_\R\{\re (\omega),\im (\omega)\} \subseteq H^1(X_\omega;\R)$. These results help us rule out \veigenvalues~for almost all directions before any splitting occurs. It is important to highlight that, although the results in \cite{CF} represent a significant novelty with respect to~\cite{AvFo} and the subsequent works~\cite{AvLe,AFS23}, in which the splitting is controlled by a probabilistic exclusion argument, the paper \cite{CF} deals essentially with a single \veigenvalue~and not with an interval of \veigenvalues. 

Nevertheless, since the projection of any interval of vectors $\alpha \im (\omega) - \vec{n} \in H^1(X;\mathbb{R})$ for $\vec{n} \in H^1(X;\mathbb{Z})$ onto the symplectic complement of the tautological plane, which is equal to the projection of only $- \vec{n}$, is a parallel section, the arguments of \cite{CF} can be adapted, as long
as there is no  splitting of the interval of  \veigenvalues~ caused by expansion in the unstable direction of the tautological plane. New ideas need to be introduced to deal with this splitting issue.

To prevent splitting we prove via explicit geometric arguments, motivated by work of D.~Davis and the second author \cite{CD}, that there exists a sequence of \emph{partial rigidity times} $0 <t_1<t_2<...$ with bounded ratio, i.e.,
\[
\underset{j \to \infty}{\limsup}\, \frac{t_{j+1}}{t_j}<\infty.
\]
From this we derive additional restrictions on \veigenvalues~  which imply that, for almost every direction, it is enough to consider only one of the subintervals that arises from splitting, while, at the same time, the projection onto the symplectic complement of the tautological space grows.

The arguments outlined above rule out all eigenvalues for almost every direction on any given translation surface but clearly fails for $(X, \omega)$ if there exist $\alpha\in \R\setminus \{0\}$ and $\theta \in [0,2\pi)$ such that $\alpha \im (e^{i\theta} \omega) \in H^1(X, \Z)$, i.e., if the tautological plane of $(X,\omega)$ contains a non-zero integral cohomology class.
Conversely, when there is a non-zero integer vector in the tautological plane, we build an affine map to the circle, which gives a non-trivial continuous eigenfunction in all but one direction, and a non-constant invariant function in the remaining direction.

\section{Acknowledgments}

JC acknowledges support of NSF grants DMS-2055354 and DMS-2350393.
GF acknowledges support of NSF grant DMS-2154208. We thank Pedram Safaee for a helpful, clarifying comment on an earlier draft.

\section{Statement of the main theorem and reductions}

\subsection{Outline of this section.}

After recalling in section \ref{subsec:background} background on the Hodge bundle over the moduli space of translation surfaces and on the Kontsevich-Zorich renormalization cocycle, and establishing notation, we give in section~\ref{subsec:main_th} a complete statement of our main theorem  characterizing translation surfaces whose translation flows are weakly mixing in almost all directions. 
We then introduce in section~\ref{subsec:criterion} a criterion for weak mixing which upgrades the Veech criterion with a condition derived from partial rigidity of
translation flows, based on the notion of rigidity configurations.
Finally, in sections \ref{subsec:red_lyap} and \ref{subsec:red_horo} we make two reductions. We show that for our problem it is sufficient to consider
the restriction of the Kontsevich--Zorich cocycle to a hyperbolic rational subbundle given by the sum of Galois conjugates of the tautological subbundle.  
We also make the straightforward reduction that when eliminating ``bad'' (non weak mixing) parameters, it is equivalent to consider circle orbits and
horocycle orbits of a translation surface.

\subsection{Background}
\label{subsec:background}
Let $\mathcal{M}$ be a moduli space of Riemann surfaces. Given $X \in \mathcal{M}$, consider the decomposition of its complex absolute cohomology into holomorphic and antiholomorphic 1-forms, that is,
\[
H^1(X;\C) = H^{1,0}(X) \oplus H^{0,1}(X).
\]
Denote by $\star \colon H^1(X;\C) \to H^1(X;\C)$ the Hodge star operator, i.e., the unique complex  linear operator on complex absolute cohomology that acts as $(-i)$ on holomorphic $1$-forms and as $i$ on antiholomorphic $1$-forms. Denote by $\cup$ the cup product on $H^1(X;\C)$; this is a complex symplectic form. Given $\alpha,\beta \in H^1(X;\C)$ their Hodge inner product is given by
\[
\left\langle \alpha,\beta \right\rangle_X = \alpha \cup \overline{\star \beta};
\]
this is a Hermitian inner product. The (symplectic) orthogonal of a subspace $V \subseteq H^1(X;\C)$ will be denoted by $V^\perp$. The Hodge norm $\|\cdot \|_X$ on $H^1(X;\C)$ is the norm induced by the Hodge inner product $\left\langle \cdot, \cdot\right\rangle_X$.
Denote by $d_X$ the metric induced by the Hodge norm $\|\cdot\|_X$ on $H^1(X;\R)$. 

\medskip

\noindent 
{\bf Notation}: In the following we let for simplicity $d$ denote the metric induced by the Hodge norm $\Vert \cdot \Vert$ on the real cohomology bundle $H^1_\R$ with fiber $H^1(X,\R)$. Namely,
for any $X\in \mathcal{M}$ and any pair of vector $\vec{w}, \vec{w} \in H^1(X, \R)$ we denote
$$
d(\vec{v}, \vec{w}) = \Vert \vec{v}- \vec{w} \Vert := \Vert \vec{v}- \vec{w} \Vert_X\,.
$$
Abelian differentials considered will always be holomorphic. Let $\mathcal{H}$ be a stratum of Abelian differentials. 
Recall that the linear action of $2 \times 2$ matrices on the plane induces a $\mathrm{GL}^+(2,\R)$-action on $\mathcal{H}$. This action restricts to an $\mathrm{SL}(2,\R)$-action on $\mathcal{H}^1 \subseteq \mathcal{H}$, the corresponding locus of unit area Abelian differentials. For $t \in \R$, $s \in \R$, and $\theta \in [0,2\pi]$, denote
\[
g_t := \left(\begin{array}{cc}
   e^t  & 0 \\
    0 & e^{-t}
\end{array} \right), \quad 
h_s := \left( \begin{array}{cc}
    1 & s \\
    0 & 1
\end{array}\right), \quad
r_\theta := \left( \begin{array}{cc}
    \cos \theta & -\sin \theta \\
    \sin \theta & \cos \theta
\end{array}\right).
\]

The bundles over $\mathcal{H}^1$ whose fibers above each $(X,\omega) \in \mathcal{H}^1$ are given by $H^1(X;\C)$ and $H^1(X;\R)$ are denoted by $H^1_\C$ and $H^1_\R$ respectively; unless otherwise stated, the fibers of these bundles will always be endowed with the corresponding Hodge inner products, Hodge norms, and cup products. 

Recall that the $\mathrm{SL}(2,\R)$-action on $\mathcal{H}^1$ extends to $H^1_\C$ and {$H^1_\R$} by parallel transport with respect to the (flat) Gauss-Manin connection; unless otherwise stated, when refering to parallel transport, this will always be done with respect to this connection. The action of the subgroup $\mathcal{G} := \{g_t\}_{t \in \R}$ on $H^1_\R$ or $H^1_\C$ is the (complexified) Kontsevich-Zorich  cocycle. Unless otherwise stated, when referring to the Lyapunov exponents of an $\mathrm{SL}(2,\R)$-invariant sub-bundle of $H^1_\R$ or $H^1_\C$, this will always be done so with respect to the Kontsevich-Zorich cocycle and some appropriate $\mathrm{SL}(2,\R)$-invariant ergodic measure on the base $\mathcal{H}^1$.

Given a translation surface $(X,\omega) \in \mathcal{H}$, denote by $|\omega|$ the measure induced by $\omega$ on $X$, {which is given by the smooth non-negative $2$-form
$$
\frac{i}{2} \, \omega \wedge \overline{\omega} =\re \omega\wedge \im \omega \,,
$$
vanishing precisely at the zero set $\Sigma$ of $\omega$.} Recall that, when convenient, we record the position of the zero set $\Sigma \subseteq X$ of an Abelian differential $\omega$ on $X$ by denoting it by $(X,\Sigma,\omega) \in \mathcal{H}$.
For any translation surface, that is a pair $(X,\omega) \in \mathcal{H}^1$, let
$L^2(X, \omega)$ denote the Hilbert space of square-integrable
complex values functions integrable with respect to the area measure $|\omega|$. 

\subsection{Statement of the main theorem} 
\label{subsec:main_th}
Given a translation surface $(X,\omega)$, its tautological plane is the plane in real absolute first cohomology spanned by $\re (\omega)$ and $\im (\omega)$, i.e.,
\begin{equation}
\label{eq:taut}
\mathrm{Taut}(X,\omega):= \mathrm{span}_\R\{\re (\omega),\im (\omega)\} \subseteq H^1(X;\R).
\end{equation}
Unless otherwise stated, when referring to weak mixing of a translation flow on a translation surface, this will be done so with respect to the natural measure induced by the corresponding Abelian differential. Furthermore, when referring to almost every direction on a translation surface, this will be done so with respect to the Lebesgue measure. The following is the main result of this paper.

\begin{thm}
    \label{thm:main1}
    Let $(X,\omega)$ be an arbitrary translation surface. Then, the following statements are equivalent:
    \begin{enumerate}
        \item The translation flow on $(X,\omega)$ is weak mixing in almost every direction.
        \item The tautological plane of $(X,\omega)$ contains no non-zero integer points, i.e., 
        \[
        \mathrm{Taut}(X,\omega) \cap H^1(X;\Z) = \{0\}.
        \]
        \item The translation surface $(X,\omega)$ admits no affine map to the circle $\mathbb{S}^1 = \R/\Z$ endowed with its natural translation structure.
        \item The translation flow on $(X,\omega)$ is weak mixing in some direction.
        \item The translation flow on $(X,\omega)$ is weak mixing in a $G_\delta$, i.e. countable intersection of open sets, dense set of directions.
    \end{enumerate}
\end{thm}

Most of the implications in Theorem \ref{thm:main1} follow via standard arguments, with the exception of (2) implies (1), which encompasses the bulk of the proof. The implication $(5) \Rightarrow (4)$ is immediate and  $(1) \Rightarrow (5)$ follows from the standard fact that weak mixing is a $G_\delta$ Property. In addition, we have:

\begin{proof}
$(4)\Rightarrow (3):$ We prove the contrapositive statement that if the translation surface admits an affine map to the circle $\T = \R/ \Z$, then there is no direction such that translation flow is weakly mixing.  In fact, any translation 
flow on $(X, \omega)$ projects under the affine map to a linear
flow on $\T$, which is either constant (generated by the zero vector field), or non-constant. In the first case, the translation flow on  $(X, \omega)$ leaves invariant the fibers of the affine map to $\T$, which are generically circles, hence it is not ergodic, hence not weakly mixing. In the second case, the pull-back of any non-constant eigenfunction
for the linear flow on $\T$ is a non-constant eigenfunction for
the translation flow on $(X, \omega)$, hence the latter is not
weakly mixing.
\end{proof}

\begin{proof}
$(3)\Rightarrow (2)$:
We prove the contrapositive statement that if the tautological plane of $(X,\omega)$ contains a non-zero integer point, then
the translation surface $(X, \omega)$ admits an affine map to
the circle $\T = \R/ \Z $.

Let us assume that there exists $a, b\in \R$ such that
$[a \re (\omega) + b \im (\omega)] \in H^1(X, \Z)$. Let
$\Hol_\omega: X \to \R / \Z$ denote the map defined as follows. Let $x_0 \in X$ be a fixed (regular) point and let
$$
\Hol_\omega (x) = \int_{x_0}^x  a \re (\omega) + b \im (\omega) \,, \quad \text{ for all } x \in X\,.
$$
The above map is well-defined independently of the choice of
an oriented path with endpoints $x_0$ and $x$ (with boundary 
$x-x_0$). In fact, if $\gamma_1$ and $\gamma_2$ are oriented paths both with boundary $x-x_0$, then $\gamma_2-\gamma_1$ is
a cycle, so that, by the assumption that $[a \re (\omega) + b \im (\omega)] \in H^1(X, \Z)$,
$$
\int_{\gamma_2-\gamma_1} a\re (\omega) + b\im (\omega) \in \Z\,.
$$
The map $\Phi$ is by definition an affine map from $(X, \omega)$ to $(\T, d\theta)$ with differential 
\[
d\Phi  = a \re (\omega) + b \im (\omega)\,. \qedhere
\]
\end{proof}

Most of this paper is devoted to the proof of the crucial implication $(2)\Rightarrow (1)$


\subsection{Criterion for weak mixing}
\label{subsec:criterion}
Our criterion for weak mixing requires a preliminary definition.

\begin{figure}
    \centering
    \includegraphics{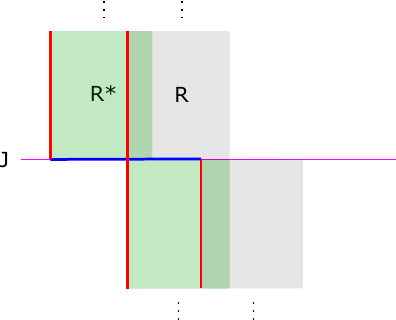}
    \caption{Rigidity configuration on a translation surface.}
    \label{fig:rigidity}
    \end{figure}

\begin{defin}
    Let $(X,\omega)$ be an arbitrary translation surface. By an immersed flat rectangle $R \hookrightarrow (X,\omega)$ we mean a smooth immersion of an open horizontal-vertical Euclidean rectangle $R \subseteq \R^2$ into $X$ which is a local isometry with respect to the singular Euclidean metric induced by $\omega$ on $X$ and whose image contains no singularities of $\omega$. Given parameters $V,H,\sigma,L > 0$, a \emph {$(V,H,\sigma,L)$-rigidity-configuration} of $(X,\omega)$ is a pair $(J,R)$, where $J \subseteq X$ is a segment of a horizontal leaf of $(X,\omega)$ and $R \hookrightarrow (X,\omega)$ is an immersed flat rectangle satisfying the following conditions:
\begin{enumerate}
    \item $J$ can be flowed upwards up to time $L$ without hitting singularities of $\omega$.
    \item Each horizontal side of $R$ maps bijectively to a subsegment of $J$.
    \item The vertical sides of $R$ have length $V$.
    \item  The distance along $J$ of the top and bottom points of every maximal vertical segment of $R$ is $H$.
    \item The maximal embedded flat sub-rectangle $R^* \subseteq R$ starting from the base of $R$ has area $|\omega|(R^*) = \sigma$.
\end{enumerate}
\end{defin}
See Figure \ref{fig:rigidity} for an example. Every rigidity configuration $(J,R)$ on a translation surface $(X,\omega)$ gives rise to a so-called `rigidity curve' by joining the top and bottom points of any maximal vertical segment of $R$ along $J$; in Figure \ref{fig:rigidity} this corresponds to the curve obtained by joining the red and blue segments.

\smallskip
\noindent {\bf Notation}:  Given any translation surface 
$(X,\omega)$ we will adopt the following notation: for all $t\in \R$
$$
(X_t, \omega_t) := g_t (X,\omega) \,.
$$



\begin{defin}\label{def:eval crit bis} Let $(X, \omega) \in \mathcal {H}$ be a translation surface. Given $\epsilon>0$ and $\mathcal{K} \subset \mathcal{H}$, we say that $\alpha\in \R\setminus \{0\}$ is an $(\epsilon, \mathcal{K},t)$-\veigenvalue~on $(X, \omega)$ 
    if both of the following two conditions hold
\begin{enumerate} 
\item\label{cond:rigid} for all $(V,H,\sigma,L)$-rigidity configuration such that $V\geq e^t$ at least one of the following hold:  $$ 
\sigma < \epsilon \quad\text{or} \quad d_\R(\alpha V ,\Z)<\epsilon \,;$$
\item\label{cond:veech} (Veech Criterion) for all $L>t$ such that 
$(X_L,\omega_L)  \in \mathcal{K}$,  $$d\Big(g_L \cdot \alpha \im (\omega), H^1(X_L;\Z) \Big) < \epsilon;$$
\end{enumerate}
\end{defin}

For readability we will often use $H^1(\Z)$ for 
$H^1(X;\Z)$  when the meaning is clear.


In the next lemma and throughout the remainder of the paper, the word {\it eigenvalue} means an eigenvalue for a {\it measurable} eigenfunction.
\begin{lem}
\label{lem:candidate pos}
If the vertical flow is ergodic and $\alpha$ is an eigenvalue for the vertical flow on $(X,\omega)$, then for every $\epsilon>0$ and compact $\mathcal{K} \subset \mathcal{H}^1$ there exists $t_0:= t_0(\epsilon, \mathcal{K}, \alpha)>0$ such that $\alpha$ is a $(\epsilon,\mathcal{K},t)$-\veigenvalue~for all $t > t_0$. 
\end{lem}

\begin{proof} 
It suffices prove the following two statements:  

\begin{itemize}
\item[a)] If there exists $\epsilon>0$, $\mathcal{K}$ compact and a diverging sequence $(t_i)\subset \R$,  such that for all $i\in \mathbb{N}$ there is an $(e^{t_i},H,\sigma,L)$-rigidity configuration with $\sigma\geq \epsilon$ and $$d_\R(\alpha e^{t_i} ,\Z)>\epsilon\,,$$  then $\alpha$ is not an eigenvalue.  

\item[b)] If there exists $\epsilon>0$, $\mathcal{K}$ compact and a diverging sequence $(t_i)\subset \R$,  such that for all $i\in \mathbb{N}$ we have  $(X_{t_i},\omega_{t_i}) \in \mathcal{K}$  
  and $$d\Big(g_{t_i} \cdot \alpha \im (\omega), H^1(X_{t_i};\Z) \Big)  > \epsilon\,,$$ then $\alpha$ is not an eigenvalue.
\end{itemize}
Indeed, because in Definition \ref{def:eval crit bis} there are only two cases, if there exists $\epsilon>0$ and $\mathcal{K}$ compact such that, for all $t_0 >0$ there exists $t>t_0$ with the property that $\alpha$ is not a $(\epsilon,\mathcal{K},t)$ \veigenvalue~for all $t>0$,  then there is a diverging sequence $(t_i)$ as in a) or such a sequence as in b). 

\smallskip
We first treat a): observe that if $(J,R)$ is a $(V,H, \sigma,L)$-rigidity configuration and $x \in R$ then $F^Vx$ is on a horizontal segment of containing $x$ and its displacement along this segment is $H$. Thus there is a set of measure $\sigma>0$, $R^*$, so that $d(F^Vx,x)\leq |H|$ for all $x \in R^*$. By the minimality of the vertical flow, there is a dense vertical separatrix and so for every $\delta>0$ there exists a $V_0$ so that if $V\geq V_0$ then the horizontal component of an immersed Euclidean rectangle with vertical component $V$ is less than $\delta$. 
Now let $f$ be a measurable eigenfunction with eigenvalue $\alpha$ and let $K\subset X$ be a compact continuity set of $f$ with measure $\vert \omega\vert(K) \geq 1-\frac{\sigma}4$ (which exists by Lusin's theorem). Let $\delta>0$ be so that if $x,y \in K$ and $d(x,y)<\delta$ then $|f(x)-f(y)|<\frac \epsilon 4$. Thus if $x \in R^*$, $x,F^Vx$ are in $K$ and $\vert H \vert<\delta$ (since by ergodicity we can assume $\vert f \vert=1$ almost everywhere),
$$|f(F^V x)-f(x)|=|e^{iV\alpha}-1|<\frac{\epsilon}4\,.$$ 
By our choice on the size of $K$, a set of measure at least $\sigma-2\frac{\sigma}4>0$ satisfies the above conditions thereby giving the claim.

\smallskip
We then treat b): in this case the  statement is given by \cite[Theorem 4.12]{FCet} unless the eigenfunction is actually continuous. In this case there is no need to use Lusin sets in the proof of the Veech criterion, which therefore holds for any given compact subset $\mathcal{K} \subset \mathcal{H}^1$.
\end{proof}

\begin{lem}\label{lem:candidate mble} For every $J \subset [-1,1]$ measurable with $\psi:J \to \R\setminus \{0\}$ measurable, for  all $t \in \R^+,\epsilon>0$ and $\mathcal{K} \subset \mathcal{H}^1$, compact, the set $s \in J$ so that $\psi(s)$ is a  $(\epsilon, \mathcal K, t)$-\veigenvalue~for the vertical flow on $h_s(X,\omega)$ is measurable.
\end{lem}

\begin{proof} 
Let (1) and (2) be the properties
in the definition of a \veigenvalue~(see Def. \ref{def:eval crit bis}). It suffices to show that the sets 
$$
\begin{aligned} 
(i) &\quad \{s \in J:\psi(s) \text{ satisfies (1) for some } t'\geq t\}  \,, \\
(ii) &\quad \{s \in J:\psi(s) \text{ satisfies (2) for some }t'\geq t\}
\end{aligned}$$ 
are both measurable. The first set is measurable because the set 
$$
\{ (V,\alpha) \in \R^+ \times \R\setminus \{0\} \vert d_\R (V \alpha,\Z) <\epsilon\}
$$
is open and the set of all $V>e^t$ such that there exists a
$(V,H,\sigma,L)$ rigidity configuration with $\sigma > \epsilon$ is countable and depends continuously on $s\in J$. In fact, for any rigidity configuration
$V$ is the vertical length of a transversal closed curve, homotopic
to the union of a vertical segment of length $V$ and a horizontal segment of length $H$, hence it is a period of the imaginary part of the Abelian differential. The sets of such periods is countable since the set of homotopy classes of loops is a countable set. Each of such periods depends continuously on the parameter $s\in J$. Moreover,  for any fixed homotopy class, its vertical period and the maximal area of a rectangle $R^*$ of the associated rigidity configuration vary continuously with the parameter $s\in J$. In conclusion, 
there exists a sequence of measurable functions $(V_n)$  defined on $J$ such that $V_n$ is the vertical length of a rigidity configuration on the open set $J_n:=\{s\in J\vert V_n(s) \not =0\}$. Such functions are given by vertical periods on $h_s(X, \omega)$ of loops, extended by the value $0$ outside of the interval on which the vertical period does not come from a rigidity configuration. 

Then we have that 
 $$
\begin{aligned}
\{s \in J:&\psi(s) \text{ satisfies (1) for some } t'\geq t\} \\ & \quad = \bigcap_{n\in \N}  (V_n\times \psi)^{-1} \{ (V,\alpha) \in  \R^+ \times \R\setminus\{0\} \vert V > e^t \text{ and } d_\R(V \alpha, \Z) < \epsilon\}\,, 
\end{aligned}
$$
hence the set in $(i)$ is measurable.


\smallskip 
\smallskip 
For the second set for each $q\in \Q^+$ and $n\in \mathbb{N}$ let 
$$W_{q,n}(\epsilon)=\{(X, \omega, \alpha):d(g_q \cdot \alpha \im (\omega), H^1(X_q;\Z))> \epsilon\text { and } g_q(X, \omega) \in B(\mathcal{K},
\frac 1 n)\}.$$ 
This is an open set and the set we are considering is the
set of $s \in J$ such that
\[
(h_s(X, \omega), \psi(s) ) \in 
\bigg(\bigcup_{r=0}^\infty 
\bigcap_{k=1}^\infty \bigcap_{n=1}^\infty 
\bigcup_{q\in [t+r,t+r+1] \cap\Q} W_{q,n}(\epsilon-\frac 1 k)\bigg)^c. \qedhere
\]

\end{proof}

\subsection{Reduction to positive Lyapunov exponents}\label{subsec:red_lyap} 
In working with Definition \ref{def:eval crit bis} \eqref{cond:veech} it is enough to work with $\mathrm{SL} (2,\R)$-invariant, strongly irreducible subbundles of $H^1_\R$ whose top Lyapunov exponent is positive; this guarantees expanding properties of the Kontsevich-Zorich cocycle in generic directions, leading to a contradiction with the statement in \eqref{cond:veech}.  {The reduction to this non-uniformly (partially) hyperbolic situation follows from} 
known facts about the Galois theory of affine invariant sub-manifolds, which we now recall.

Recall that an affine  $\mathrm{GL}(2,\R)$-invariant submanifold $\mathcal{N} \subseteq \mathcal{H}$ is an immersed suborbifold of $\mathcal{H}$ locally defined in period coordinates by homogeneous equations with real coefficients. Denote by $k(\mathcal{N}) \subseteq \R$ the smallest real subfield containing the coefficients of any such set of equations, i.e., the field of definition of $\mathcal{N}$; this definition makes sense because the lattice $H^1(X,\Sigma;\Z)$ corresponds to integer points in period coordinates. In general, the field of definition of a subspace/subbundle of a real/complex vector space/bundle with a prescribed integer lattice is the smallest real/complex subfield over which homogeneous equations cutting out the subspace/subbundle can be found.
It is equal to the intersection of the holonomy fields of all translations surfaces in $\mathcal{N}$ \cite[Theorem 1.1]{FOD}.

Denote by $\mathcal{N}^1 \subseteq \mathcal{H}^1$ the unit area locus of $\mathcal{N}$. Recall that, via period coordinates, the tangent space of $\mathcal{N}$ at any translation surface $(X,\Sigma,\omega) \in \mathcal{N}$ is naturally identified with a subspace\footnote{Note that, while $\mathrm{SL}(2,\mathbb{R})$ orbit closures need not be manifolds, because self-intersection loci and orbifold loci are closed and $\mathrm{SL}(2,\mathbb{R})$-invariant, $\omega$ is a regular point of $\overline{\mathrm{SL}(2,\mathbb{R})\omega}$.} 
\[
T\mathcal{N}_{(X,\Sigma,\omega)} \subseteq H^1(X,\Sigma;\C).
\] 
For convenience we consider the bundles
\[
H^1_{\mathrm{rel},\C}, \quad T\mathcal{N} \subseteq H^1_{\mathrm{rel},\C}|_\mathcal{N}
\]
over $\mathcal{H}^1$ and $\mathcal{N}^1$ whose fibers above every $(X,\Sigma,\omega) \in \mathcal{H}^1$ are given by
\[
H^1(X,\Sigma;\C), \quad T\mathcal{N}_{(X,\Sigma,\omega)} \subseteq H^1(X,\Sigma;\C),
\]
respectively. Notice that, because $\mathcal{N}$ is an affine invariant submanifold, $T\mathcal{N}$ is invariant under parallel transport. Notice also that, essentially by definition, the field of definition of $T\mathcal{N}$ is $k(\mathcal{M})$. Denote by $p \colon H^1_{\rm{rel},\C} \to H^1_\C$ the natural map which on fibers send relative cohomology to absolute cohomology.

By the work of Eskin, Mirzakhani, and Mohammadi \cite{EMM, EM}, every $\mathrm{GL}^+(2,\R)$ orbit closure in $\mathcal{H}$ is an affine invariant submanifold $\mathcal{N}$ whose unit area locus $\mathcal{N}^1$ supports a natural $\mathrm{SL}(2,\R)$-invariant ergodic probability measure $\mu$; in an abuse of terminology we refer to this measure as the affine measure of $\mathcal{N}^1$. Given a translation surface $(X,\omega) \in \mathcal{H}^1$ whose translation flow (in some direction) we would like to study, it is natural to consider its orbit closure and the corresponding affine measure. Unless otherwise stated, when referring to the Lyapunov exponents of the Kontsevich-Zorich cocycle of an $\mathrm{SL}(2,\R)$-invariant subbundle of $H^1_\R$ or $H^1_\C$ over an $\mathrm{SL}(2,\R)$-orbit-closure, we will always consider the corresponding affine measure as the underlying measure.

The field of definition of a flat subbundle $E \subset  
H^1_{\C}$ is the smallest subfield of $\R$ so that locally the linear subspace $E$ of $H^1(X, \C)$
can be defined by linear equations (with respect to an integer basis of
$H^1(X, \Z)$)  with coefficients in this field. The trace field of a flat bundle
over $\mathcal M$ is defined as the field generated by traces of the corresponding
representation of the fundamental group $\pi_1(\mathcal N)$.

Given a subspace/subbundle $V$ of a complex vector space/bundle $W$ with a prescribed integer lattice whose field of definition is $k \subseteq \C$, it is natural to consider the subspace/subbundle in $W$ spanned by all the Galois conjugates of the subspace $V$, i.e., by all the subspaces obtained by applying the natural extension of every field embedding $k \hookrightarrow \C$ to $V \subseteq W$. In this context, we consider the following theorem of A.~Wright.

\begin{thm} \cite[Theorem 1.5]{FOD}
    \label{thm:galois_wright}
    Let $\mathcal{H}$ be a stratum of Abelian differentials and $\mathcal{N} \subseteq \mathcal{H}$ be an affine invariant submanifold with field of definition $k(\mathcal{N}) \subseteq \R$. Then, the field of definition of $p(T\mathcal{N}) \subseteq H^1_\C|_{\mathcal{N}^1}$ is also $k(\mathcal{N})$. Also, if $\Gal(\mathcal{N}) \subseteq H^1_\C|_{\mathcal{N}^1}$ denotes the subbundle spanned fiberwise by all the Galois conjugates of $p(T\mathcal{N})$, then the field of definition of $\Gal(\mathcal{N})$ is $\Q$. 
\end{thm}

As Theorem \ref{thm:galois_wright} guarantees $\Gal(\mathcal{N})$ is defined over $\Q$, there exists a subbundle $\Gal(\mathcal{N})_\R \subseteq H^1_\R|_{\mathcal{N}^1}$ with field of definition $\Q$ such that
\[
\Gal(\mathcal{N}) = \Gal(\mathcal{N})_\R \otimes_\R \C.
\]
Notice that $\Gal(\mathcal{N})$ and $\Gal(\mathcal{N})_\R$ are invariant with respect to parallel transport.

\smallskip
Let $\Lambda(\mathcal{N}) \subseteq H^1_\R|_{\mathcal{N}^1}$ the fiberwise lattice over ${\mathcal{N}^1}$ whose fiber above every translation surface $(X,\omega) \in \mathcal{N}^1$ is given by the intersection of the fiber of $\Gal(\mathcal{N})_\R$ at $(X, \omega)$ with the integer lattice $H^1(X;\Z) \subseteq H^1(X;\R)$. This construction gives rise to a fiberwise lattice precisely because the field of definition of $\Gal(\mathcal{N})_\R$ is $\Q$. Furthermore, this fiberwise lattice is invariant with respect to parallel transport.
\begin{lem} 
\label{lem:lattices}
For any compact set $\mathcal{K} \subset \mathcal{H}$ there
exists a constant $c_\mathcal{K}>0$ such that, for every $X  \in \mathcal{K}$ and every $\vec{v} \in \Gal(\mathcal{N})_\R \cap H^1(X, \R)$ we have
$$
d \Big(\vec{v}, H^1_\Z\Big) \geq c_\mathcal{K} d \Big(\vec{v}, \Lambda(\mathcal{N}) \Big)\,.
$$
\end{lem}
\begin{proof} Let $H^1_\T$ denote the toral bundle
over the moduli space $\mathcal {M}$ of Riemann surfaces,
that is, the quotient bundle $H^1_\T:= H^1_\R/ H^1_\Z$. The projection of the subbundle
$\Gal(\mathcal{N})_\R$ onto $H^1_\T$
is a closed (toral) subbundle since $\Gal(\mathcal{N})_\R$ is defined over $\Q$. 
Let us assume that the inequality in the statement does not hold, then there exists a compact set $\mathcal K\subset \mathcal H$ and a sequence $(X_n, \vec{v}_n)$ such that $X_n \in  \mathcal{K}$ and $\vec{v}_n\in \Gal (\mathcal{N})_\R \cap H^1(X_n, \R)$ such that 
$$
d \Big(\vec{v}_n, H^1(X_n, \Z) \Big) \leq  d \Big(\vec{v}_n, \Lambda(\mathcal{N}) \cap H^1(X_n, \R)  \Big) /n \,.
$$
Let $\pi_\T:H^1_ \R \to H^1_\mathbb{T} $ denote the canonical projection.  Since the restriction $H^1_\mathbb{T} \vert \mathcal{K}$ is a compact space, there exist $(X,v)$
with $X\in \mathcal{K}$ and $\vec{v} \in \Gal (\mathcal{N})_\R \cap H^1(X, \R)$ such that
after passing to a subsequence
$$
X_n \to X \,, \quad \pi_\T(\vec{v}_n) \to  \pi_\T(v) \in H^1(X, \mathbb{T})  \,.
$$
By the above inequality we have that
$\pi_\T(\vec{v}_n) \to 0$, hence $\pi_\T(v)= 0 \in H^1(X, \mathbb{T})$.
Since $0 \in \pi_\T (\Lambda(\mathcal{N}))$, it follows that
$$
d_\mathbb{T} \Big(\pi_\T(\vec{v}_n), \pi_\T (\Lambda(\mathcal{N}) )\Big) \to 0 
$$
and, moreover, since $\pi_\T (\Lambda(\mathcal{N}) ) \subset
\pi_\T (H^1_\Z ) =\{{\bf 0}\}$, we have that
$$
d_\mathbb{T} \Big(\pi_\T(\vec{v}_n), \pi_\T (\Lambda(\mathcal{N}) )\Big) = d_\mathbb{T} \big(\pi_\T(\vec{v}_n), {\bf 0}\big) 
$$
The latter identity is equivalent to
$$
d \Big(\vec{v}_n, \Lambda(\mathcal{N}) \cap H^1(X_n, \R)  \Big)=d \Big(\vec{v}_n, H^1(X_n, \Z) \Big)\,,
$$
hence we have reached a contradiction and the statement is proved.
\end{proof}

By Lemma~\ref{lem:lattices}, since for every translation surface $(X, \omega) \in \mathcal{N}$, the line $\R \im(\omega) \subset \Taut(X,\omega) \subset \Gal(\mathcal{N})$, it follows that whenever (by the Veech criterion)
$$
\lim_{\substack{t \to \infty, \\ X_t \in \mathcal{K}}} d\Big(g_t \cdot \alpha \im \omega, H^1(X_t;\Z)\Big) = 0\,,
$$
we also have
$$
\lim_{\substack{t \to \infty, \\ X_t \in \mathcal{K}}} d\Big(g_t \cdot \alpha \im \omega, \Lambda\Big) = 0 \,.
$$
With this in mind, 
we state the following result of S.~Filip, which states that the subbundle $\Gal(\mathcal{N}) \subseteq H^1_\C|_{\mathcal{N}^1}$ has no zero Lyapunov exponents; this guarantees expanding properties for the Kontsevich-Zorich cocycle in generic directions, leading to a contradiction with the Veech criterion.

\begin{thm} \cite[Corollary 1.3]{mono}
    \label{thm:galois_filip}
    Let $\mathcal{H}$ be a stratum of Abelian differentials, $\mathcal{N} \subseteq \mathcal{H}$ be a $\mathrm{GL}^+(2,\R)$-orbit-closure, and $\mu$ be the corresponding affine measure on $\mathcal{N}^1$. Then, the Lyapunov exponents of the Kontsevich-Zorich cocycle on the subbundle $\Gal(\mathcal{N}) \subseteq H^1_\C|_{\mathcal{N}^1}$ with respect to the affine measure $\mu$ are all non-zero.
\end{thm}

\subsection{From rotations to horocycles}
\label{subsec:red_horo}

\begin{prop}
\label{prop:selection}
    Let $\mathcal{H}$ be a stratum of Abelian differentials over a moduli space of Riemann surfaces $\mathcal{M}$ and $(X,\omega') \in \mathcal{H}$ be a translation surface whose translation flow in a set of directions of positive Lebesgue measure is not weak mixing. Then, there exists $\psi \in [0,2\pi)$, 
    such that, for $(X, \omega) = r_{\psi} (X, \omega')$, the set $\mathcal{S} \subseteq [-1,1]$ of $s \in [-1,1]$ for which the vertical flow $\Phi(s) :=\{\phi_t(s)\}_{t\in \R}$ of $h_{s}(X,\omega) \in \mathcal{H}$ is ergodic, but not weak mixing, has positive Lebesgue measure. Furthermore, there exists a measurable function 
    \[
    \alpha \colon {\mathcal{S}} \to \R \setminus \{0\}
    \]
    such that for every $s \in \mathcal{S}$ there is  $f(s) \in L^2(X, \omega)$, an  
    eigenfunction of the vertical flow of $h_{s} (X, \omega) \in \mathcal{H}$ with eigenvalue $\alpha(s) \neq 0$, i.e., for every $t \in \R$, the following holds,
    \[
    f(s) \circ \phi_t(s) = e^{2\pi i \alpha(s)t} \cdot f(s).
    \]
\end{prop}
To prove Proposition \ref{prop:selection} we establish in Corollary \ref{cor:generator eval mble} that for a 
one-parameter strongly continuous family of unitary flows, the pairs of the parameter and the set of eigenvalues of its infinitesimal generator is Borel. We were surprised that we did not find this in the literature, though that may be from not knowing the correct references. If this is not known, perhaps the reason is that for applications, there are general arguments that are good enough. For example, it is easy to show the set is analytic and for our purposes we could have used the Jankov-von Neumann selection theorem.


\begin{lem} 
\label{lemma:Borel}
Let $U_s$ denote the unitary (Koopman) operator on $L^2(X, \omega)$ given by the time-$1$ map of the vertical flow $\phi_1(s)$ on $h_s(X,\omega)$. The set $$\{(s,\alpha) \in [-1,1] \times \R:  e^{ i\alpha}  \text{ is an eigenvalue of } U_s \}$$ is Borel measurable. 
\end{lem}
\begin{proof} 
We use the spectral theorem.  
Denote the interval 
$\{e^{ i\theta}:\theta \in [\frac{j}{2^k},\frac{j+1}{2^k}]\}$ by $I_{j,k}$ and  let 
$g_{j,k}$ be a polynomial in $z,\bar{z}$ so that $|g_{j,k}(e^{ i\theta})-1|<\frac 1 k$ for all $\theta\in I_{j,k}$ in the interval and 
$|g_{j,k}(e^{ i\theta})|<\frac 1 k$ for all $\theta \in I_{j-1,k} \cup I_{j,k}\cup I_{j+1,k}$. Let $G_k=\{g_{j,
k}\}_{j=0}^{2^k}.$ 

Let $\{f_\ell\}$ be an orthonormal basis and $U_s$ be the unitary (Koopman) operator given by  $\phi_1(s)$. Define $$E_k(f_\ell,\epsilon)=\{\{s\}\times [\frac{j}{2^k},\frac{j+1}{2^k}]:j \in \{0,\ldots,2^{k}-1\},\, \,|\langle f_\ell,g_{{j},k} (U_{s})f_\ell\rangle|>\epsilon\}.$$ 

The set of eigenvalues of $\phi_1(s)$ is the set of atoms of the spectral measures of $U_{s}$. This is the set of pairs $(s,e^{ i \psi})$ where $\psi$ belongs to the set 
$$A_1:=\bigcup_{\ell=1}^\infty \bigcup_{n=1}^\infty \bigcup_{M=1}^\infty\bigcap_{k=M}^\infty E_k(f_\ell,\frac 1 n).$$ 
Indeed, if $\sigma_{s,f_\ell}$ is the spectral measure given by $U_s$ and $f_\ell$ then $$\bigcup_{M=1}^\infty \bigcap_{k=M}^\infty E_k(f_\ell,\frac 1 n)$$ contains any $(
s,\psi)$ so that $\sigma_{s, f_\ell}(\{e^{i\psi}\})>\frac 1 n$ and is contained in the set of $(s,\psi)$ so that $\sigma_{s,f_\ell}(\{e^{i\psi}\})\geq \frac 1 n.$

We now show that $E_k(f_\ell,\frac 1 n)$ is measurable for all $k \in \Z$ and $\ell, n \in \mathbb{N}$. 

 Since for every $h \in L^2(X, \omega)$ and for 
all $a,b\in \N$, the map $s \mapsto \langle h,U_s ^a (U_s^*)^b h\rangle= \langle h,U_s ^{a-b} h\rangle $ is continuous as a map from $[-1,1]$ to $\C$ (the Euclidean topology on $[-1,1]$ is the same as the Strong Operator Topology on the operators), and since $g_{j,k}$ is a polynomial, the set $\{s \in [-1,1] :|\langle h,g_{j,k}(U_s)h\rangle |>\epsilon\}$ is open and thus measurable. Thus its product with any interval is measurable. It follows that $E_k(f_\ell,\frac 1 n)$ is  measurable and so is $A_1$. 

\end{proof}
Let $\{U^t_{s}\}_{t\in \R}$ be the group of unitary operators on $L^2(X, \omega)$ of the flow $\{\phi_t(s)\}_{t\in \R}$ (given by composing with the elements of $\{\phi_t(s)\}_{t\in \R}$) and let $A_s$ be the self-adjoint operator such that
$iA_s$ is the infinitesimal generator of $\{U^t_{s}\}_{t\in \R}$.
\begin{cor}\label{cor:generator eval mble} The set 
$$\{(s,\alpha)\subset [-1,1]\times \R:\alpha \text{ is an eigenvalue of } A_s \}$$ is Borel.
\end{cor}
\begin{proof}
For every $s\in [-1, 1]$, let $A_s$ be the (not necessarily bounded) self-adjoint operator such that $iA_s$ is the infinitesimal generator of $\{U^t_{s}\}_{t\in \R}$. Let $(id,2^k):[-1,1] \times \R$ by $(id,2^k)(s,\alpha)=(s,2^k\alpha)$. Define 
$$E_k^{(q)}(f_\ell,\epsilon)=\{\{s\}\times [\frac{j}{2^k},\frac {j+1}{2^k}]:j \in \{0,\dots, 2^k-1\}, \, |\langle f_\ell,g_{j,k}(U_{s}^{\frac 1 {2^q}})f_\ell\rangle|>\epsilon\}.$$ The set of $(s,\psi)$ so that $\psi$ is an eigenvalue of $A_s$ is 
$$\hat{B}:=\bigcup_{\ell=1}^\infty \bigcup_{n=1}^\infty\bigcup_{N=1}^\infty\bigcap_{q=N}^\infty\bigcup_{M=1}^\infty \bigcap_{k=M}^\infty(id,2^q) \cdot E_k^{(q)}\left(f_\ell,\frac 1 n\right) $$
This claim can be proved by an argument similar to the one given above in the proof of Lemma~\ref{lemma:Borel} with in addition the following observation: 
if $\sigma_{s,f}$   is  the  spectral measure for the infinitesimal generator of $\{U_{s}^t\}_{t\in \R}$ of the function $f\in L^2(X, \omega)$, 
then $$(s,{\alpha})\in \bigcup_{M=1}^\infty \bigcap_{k=M}^\infty(id,2^q) \cdot E_k^{(q)}(f,\epsilon) \subset \left\lbrace(s,\alpha): \sigma_{s,f}\left(\bigcup_{j\in \Z}\{j2^q+\alpha\}\right)\geq \epsilon \right\rbrace.$$ 
Because $\sigma_{s,f}$ is a finite measure\footnote{
Thus, $\lim_{q \to \infty} \sigma_{s,f}\left(\bigcup_{j \in \mathbb{Z} \setminus \{0\}} \{j 2^q+\alpha\}\right) = 0$.}, $(s,\alpha)\in \hat{B}$ implies $\alpha$ is an eigenvalue for $A_s$. 
\end{proof}
\begin{proof}[Proof of Proposition \ref{prop:selection}]
Because $$\begin{pmatrix}
\cos(\theta)&-\sin(\theta)\\ \sin(\theta)&\cos(\theta)\end{pmatrix}= \begin{pmatrix}
1&0\\ \tan(\theta) & 1\end{pmatrix}\begin{pmatrix}
\sec(\theta)&0\\ 0 & \csc(\theta)\end{pmatrix}\begin{pmatrix}
1&-\tan(\theta)\\ 0 & 1\end{pmatrix},$$ for $\theta \notin \{\pm \frac{\pi}2\}$ the translation flow in direction $\theta$ on $(X,\omega)$ is not weakly mixing iff 
the vertical flow on $h_{-\tan(\theta)}(X, \omega)$ is not weakly mixing. To find $\psi$, we choose a density point of the measurable set of $\theta$ so that the flow in direction $\theta $ is not weakly mixing. So in any neighborhood of $\psi$ the set of non-weakly mixing directions has positive measure and thus for any $c>0$ we have that for a positive measure set $\mathcal S$  of $s\in (-c,c)$ the vertical flow on $h_sr_{\psi}(X,\omega)$ is not weakly mixing. 
With this in hand, by Corollary \ref{cor:generator eval mble} says that  $$B:=\{(s,\alpha)\subset [-1,1]\times \R:\alpha \text{ is an eigenvalue of  }A_s\}
$$ is a Borel set  and $\{s\}\times \R \cap B$ is at most countable. By \cite[Theorem 18.10]{Kec} 
we have a measurable section $\alpha: \mathcal{S} \to \R$ 
such that $\alpha(s)$ is an eigenvalue of $A_s$. 

\end{proof}

\section{Proof of (2) implies (1) in Theorem \ref{thm:main1}}

This section proves the  results we need for the main implication of this paper, namely that the absence of a non-zero integer vector in the tautological plane implies that the flow is weakly mixing in almost every direction. 

The first result we establish in section \ref{subsec:rigid_conf} concerns the abundance of $(V,H,\sigma,L)$-rigidity configurations. The next result in section \ref{subsec:local} interfaces the restrictions coming from rigidity configurations with the Veech criterion. In section \ref{subsec:repelling} we quote the main technical result from~\cite{CF} which proves a strong result on the growth of vectors in the so-called balanced subspace of the cohomology bundle. We also adapt this result to take advantage of our reductions in section~\ref{subsec:red_lyap}. With these ingredients in hand we complete the proof.

\subsection{The existence of $(V,H,\sigma,L)$-rigidity configurations}
\label{subsec:rigid_conf}

We now aim to prove the following proposition, which, under basic compactness assumptions, provides rigidity configurations with uniform parameters.

\begin{prop}
    \label{prop:rigidity}
    Let $\mathcal{H}^1$ be a stratum of unit area Abelian differentials. Then, for every compact set $\mathcal{K} \subseteq \mathcal{H}^1$, there exist constants $L_0 = L_0(\mathcal{K}) \geq 1 $ and $C = C(\mathcal{K}) \geq 1$, such that for every $(X,\omega) \in \mathcal{K}$ and every $L \geq L_0$, there exists a $(V,H,\sigma,L)$-rigidity-configuration $(J,R)$ of $(X,\omega)$ with the following properties:
    \begin{enumerate}
        \item $L /C \leq V \leq C \cdot L$.
        \item $\sigma \geq 1/C$.
        \item If $0 < t < \log L$ is such that $(X_t,\omega_t) \in \mathcal{K}$, then $H < C \cdot e^{-t}$.
    \end{enumerate}
\end{prop}

\begin{proof}
    Fix a compact subset $\mathcal{K} \subseteq \mathcal{H}^1$ and let $L_0 = L_0(\mathcal{K}) \geq 1$ be such that for every Abelian differential in $\mathcal{K}$, the union of all downward prongs of length $L_0$ starting from every singularity of the differential intersects every horizontal leaf; if any such downward prong first hits a singularity, we consider the corresponding vertical saddle connection instead. Such $L_0$ can be found using a covering argument over $\mathcal{K}$. From now on we consider an arbitrary parameter $L \geq L_0$.
    
    Let $(X,\omega) \in \mathcal{K}$ and denote by $\Phi_\omega := \{\phi_{\omega,t}\}_{t \in \R}$ the vertical flow of  $(X,\omega)$. We partition $X$, modulo sets of $|\omega|$-measure zero, into embedded flat open rectangles with convenient geometric and dynamical properties. Consider all downward prongs of length $2L$ starting from every singularity of $\omega$; again, if any such downward prong first hits a singularity, we consider the corresponding vertical saddle connection instead. From the start and end of every such downward prong shoot rays in all possible horizontal directions until they hit one of the prongs, potentially the same; this always happens because horizontal saddle connections decompose translation surfaces into horizontally periodic and horizontally minimal components. The connected components of the complement of the downward prongs and horizontal rays give the desired partition into embedded flat open rectangles. 
    
    Notice that, by construction, if $I$ denotes one of the horizontal sides of such a rectangle, then one never hits a singularity of $\omega$ if one flows $I$ upwards up to time $2L$; compare to the first condition in the definition of rigidity configurations. Indeed, if this was not the case, then $I$ would have to have been split by one of the vertical prongs. We warn the reader that this decomposition into rectangles is not the same as a decomposition into towers of height at least $2L$ because when flowing the base of a rectangle upwards one could meet the base of a different rectangle, but never a singularity, at a time strictly less than $2L$.

    Notice the number of rectangles is bounded by a number $c_{\mathcal{H}}$ depending only on the stratum.
    Denote by $R' \subseteq X$ the corresponding rectangle of largest area; this rectangle has area at least $1/c_\mathcal{H}$. The construction above ensures $R'$ has vertical length at most $2L$ and so it must have horizontal length at least $1/2c_\mathcal{H}L$. Denote by $J' \subseteq X$ the base of $R$ so that
    \[
    |J'| \geq 1/2c_\mathcal{H}L.
    \]

    Now let $0 < \epsilon_\mathcal{K} \leq 1$ be the minimum between $1$ and the length of the shortest saddle connections of Abelian differentials in $\mathcal{K}$. Consider the compact subset $\mathcal{K} \subseteq \mathcal{K'} \subseteq \mathcal{H}^1$ defined as 
    \[
    \mathcal{K}' := \bigcup_{0 \leq s \leq \log(2) - \log(\epsilon_\mathcal{K})} g_s \mathcal{K}.
    \]
    We now split our analysis into three different cases:
    \begin{enumerate}
        \item $(X_{\log(L)},  \omega_{\log(L)}) \in \mathcal{K'}$.
        \item $(X_{\log(L)},  \omega_{\log(L)}) \notin \mathcal{K'}$ and the minimal return time of the vertical flow $\Phi_\omega$ of $(X, \omega)$ to $J'$ is at least $L/10$.
        \item $(X_{\log(L)}, \omega_{\log(L)}) \notin \mathcal{K'}$ but the minimal return time of the vertical flow $\Phi_\omega$ of $(X, \omega)$ to $J'$ is at most $L/10$.
    \end{enumerate}

    \textit{Case 1:} $(X_{\log(L)},\omega_{\log(L)}) \in \mathcal{K'}$. As above, denote by $0 < \epsilon_{\mathcal{K}'} \leq 1$ the minimum between $1$ and the length of the shortest saddle connections of Abelian differentials in $\mathcal{K}'$. Consider the leftmost horizontal subsegment $J \subseteq J'$ of length equal to $\min\{\epsilon_{\mathcal{K}'}/2L, |J|\}$. 
    
    \noindent
{\textbf{Claim:}}  The minimal return time of the vertical flow to $J$ is at least $L\epsilon_{\mathcal{K'}}/2$. 
  \begin{proof}  Let $\gamma'$ be the simple closed curve obtained by taking a vertical trajectory that returns quickest to $J$ and then closing it up along $J$; by the Gauss-Bonnet theorem, this curve is homotopically non-trivial. 
    
    Denote by $g_{\log(L)} \gamma'$ the parallel transport of this curve to $(X_{\log(L)}, \omega_{\log(L)})$. The assumption that $(X_{\log(L)}, \omega_{\log(L)}) \in \mathcal{K'}$ guarantees the length of $g_{\log(L)} \gamma'$ is at least $\epsilon_{\mathcal{K}'}$; if not, tightening this curve to a flat geodesic would produce a saddle connection of length less than $\epsilon_{\mathcal{K}'}$. The upper bound $|J| \leq \epsilon_{\mathcal{K}'}/2L$ ensures the horizontal segment of $g_{\log(L)} \gamma'$ has length at most $\epsilon_{\mathcal{K}'}/2$. It follows that the vertical segment of $g_{\log(L)} \gamma'$ has length at least $\epsilon_{\mathcal{K}'}/2$. Pulling back to the original differential we see that the vertical segment of $\gamma'$ has length at least $L\epsilon_{\mathcal{K'}}/2$, proving the claim.
    \end{proof}

    We now construct the desired rigidity configuration. Consider the interval exchange transformation $T$ given by the first return map of the vertical flow to $J$ and denote by $I$ the longest interval of this transformation. Consider the flat rectangle $R$ obtained by flowing $I$ upwards until it returns to $J$; in this case the rectangle is embedded and not only immersed. The desired rigidity configuration is given by the pair $(J,R)$.

    We now show this pair satisfies the desired conditions for an appropriate choice of constant $C = C(\mathcal{K}) > 0$. First, by  the construction of $R'$, we know that $J$ can be flowed upwards up to time $2L$ without hitting a singularity of $\omega$. Notice also that the horizontal segment of any rigidity curve of $(J,R)$ has length at most $\vert J\vert \leq \epsilon_\mathcal{K'}/2L$. 
    The  number of intervals in the interval exchange transformation $T$ is bounded by a number $d_{\mathcal{H}}$ depending only on the stratum. In particular, 
    \[
    |I| \geq |J|/d_\mathcal{H} = \min\{\epsilon_{\mathcal{K}'}/2d_\mathcal{H}L, |J'|/d_\mathcal{H}\} \geq \min\{\epsilon_{\mathcal{K}'}/2d_\mathcal{H}L, 1/2c_\mathcal{H} d_\mathcal{H}L\}.
  \]
    Recall that the vertical length of $R$ is at least $L\epsilon_{\mathcal{K}'}/2$. Furthermore, because $\omega$ has total area $1$,  the vertical length of $R$ is at most
    \[
    1/|I| \leq \max\{2d_\mathcal{H} L /\epsilon_{\mathcal{K'}}, 2c_\mathcal{H} d_\mathcal{H} L\}.
    \]
    Finally, the total area of the rectangle $R^* = R$ can be bounded as follows,
    \[
    |\omega|(R^*) = |\omega|(R) \geq |I| \cdot L\epsilon_{\mathcal{K}'}/2 \geq \min\{\epsilon_{\mathcal{K}'}^2/4 d_\mathcal{H}, \epsilon_{\mathcal{K}'}/4c_\mathcal{H}d_\mathcal{H}\}.
    \]
    We conclude that $(J,R)$ has the desired properties if
    \[
    C = C(\mathcal{K}) \geq \max\{\epsilon_{\mathcal{K}'}/2, 2/\epsilon_{\mathcal{K}'}, 2d_\mathcal{H}/\epsilon_{\mathcal{K}'}, 2c_\mathcal{H} d_\mathcal{H}, 4d_\mathcal{H}/\epsilon_{\mathcal{K}'}^2, \epsilon_{\mathcal{K}'}/4c_\mathcal{H}d_\mathcal{H}\}.
    \]

   {\textit{Case 2:}} 
   $(X_{\log(L)}, \omega_{\log(L)}) \notin \mathcal{K'}$ and that the minimal return time of $\Phi_\omega$ to $J'$ is at least $L/10$. Furthermore, we will use the condition that $(X_t,\omega_t) \in \mathcal{K}$ for $0 < t < \log L$ introduced in the third item of the statement of the proposition. Let $J = J'$, let $T$ be the interval exchange transformation given by the first return map of the vertical flow to $J$, and let $I$ be the longest interval of this transformation. Consider the flat rectangle $R$ obtained by flowing $I$ upwards until it returns to $J$; in this case the rectangle is embedded and not only immersed. The desired rigidity configuration is given by the pair $(J,R)$. 

    Let us verify this pair has the desired properties for an appropriate choice of constant $C = C(\mathcal{K}) > 0$. First, by  the construction of $R'$, we know that $J$ can be flowed upwards up to time $2L$ without hitting a singularity of $\omega$. Now recall that $|J| = |J'| \geq 1/2c_\mathcal{H}L$ and notice that, by area considerations, $|J| \leq 10/L$. In particular, the length of the horizontal component of any rigidity curve of $(J,R)$ is at most $|J| \leq 10/L$. Furthermore,
    \[
    |I| \geq |J|/d_\mathcal{H} \geq 1/2c_\mathcal{H} d_\mathcal{H} L.
    \]
    Notice that, by our assumption on the minimal return time of $\Phi_\omega$ to $J' = J$, the vertical length of $R$ is at least $L/10$. Furthermore, area considerations ensure the vertical length of $R$ is at most
    \[
    1/|I| \leq 2c_\mathcal{H} d_\mathcal{H} L.
    \]
    Notice that the embedded flat rectangle $R$ has area at least $1/d_\mathcal{H}$ times the area of the rectangle $R'$ constructed above. In particular, the total area of the rectangle $R^* = R$ can be bounded as follows
    \[
    |\omega|(R^*) = |\omega|(R) \geq |\omega|(R')/d_\mathcal{H} = 1/{\color{red}(}c_\mathcal{H}d_\mathcal{H}{\color{red})}.
    \]
    We conclude that $(J,R)$ has the desired properties if
    \[
    C = C(\mathcal{K}) \geq \max\{10,2c_\mathcal{H}d_\mathcal{H},c_\mathcal{H}
    d_\mathcal{H},10\}.
    \]

    {\textit{Case 3:}} we construct the desired rigidity configuration under the assumptions that $(X_{\log(L)}, \omega_{\log(L)}) \notin \mathcal{K'}$ and the minimal return time of $\Phi_\omega$ to $J'$ is at most $L/10$. Recall that one never hits a singularity of $\omega$ if one flows $J'$ upwards up to time $2L$. In particular, if $0 < s < L/10$ denotes the minimal return time to $J'$, then $\phi_{\omega,s}|_{J'}$ must be a horizontal translation by some amount $w \in \R$. Furthermore, by induction, for every $n \in \mathbb{N}$ such that $n s < 2L$, $\phi_{\omega,ns}|_{J'}$ is a horizontal translation by $nw$. Denote $k := \lceil L/s \rceil - 1 \in \mathbb{N}$; notice that $ks < L$. To construct the desired rigidity configuration $(J,R)$ consider $J$ to be the union of all horizontal translations of $J'$ by $nw$ for $0 \leq n \leq k$ and $R$ to be the immersed flat rectangle obtained by flowing $J'$ upwards up to time $ks$.

    Let us verify this pair has the desired properties for an appropriate choice of constant $C = C(\mathcal{K}) > 0$. First, by construction, the horizontal segment $J$ can be flowed upwards up to time $2L-ks \geq L$ without hitting any singularities of $\omega$. Notice also that the vertical length of $R$ is
    \[
    ks =(\lceil L/s \rceil - 1)s \in [L/2,L].
    \]
    Flowing $J'$ upwards until it returns to itself yields the embedded flat rectangle $R^* \subseteq R$. In particular, the area of $R^*$ is at least  the area of $R'$, i.e.,
    \[
    |\omega|(R^*) \geq |\omega|(R') \geq 1/c_\mathcal{H}.
    \]
    To complete the proof we make use of the assumption $(X_{\log(L)}, \omega_{\log(L)}) \notin \mathcal{K'}$ to control the length of the horizontal components of rigidity curves of $(J,R)$. Notice that, under this assumption, if $L > e^t$ and $(X_t,\omega_t) \in \mathcal{K}$, then, by definition of $\mathcal{K'}$,
    \[
    L > 2 e^t/\epsilon_{\mathcal{K}}.
    \]
    We claim the horizontal component of any rigidity curve of $(J,R)$ has length
    \[
    k |w| = (\lceil L/s \rceil - 1) |w| \leq 1/L \leq \epsilon_\mathcal{K} e^{-t}/2.
    \]
    Indeed, begin by noticing that, if we assume without loss of generality that $w \geq 0$, then the subinterval of $J'$ of length $|w|$ starting from its left endpoint travels upwards injectively for time at least $L$. By area considerations,
    \begin{equation}
    \label{eq:A}
    |w| \leq 1/L \leq \epsilon_\mathcal{K} e^{-t}/2.
    \end{equation}
    In particular, recalling the definition of $k$, it follows that
    \begin{equation}
    \label{eq:B}
    k |w| = (\lceil L/s \rceil - 1)/L \leq 1/s.
    \end{equation}
    Now let $\gamma$ be a rigidity curve of $(J,R)$ and denote by $g_t \gamma$ the parallel transport of this curve to $(X_t,\omega_t) \in \mathcal{K}$. By the Gauss-Bonnet theorem, $\gamma$ is homotopically non-trivial. In particular, $g_t \gamma$ has length at least $\epsilon_\mathcal{K}$; if not, tightening this curve to a flat geodesic would produce a saddle connection of length less than $\epsilon_{\mathcal{K}}$. Furthermore, by \eqref{eq:A}, the horizontal component of $g_t \gamma$ has length at most $\epsilon_\mathcal{K}/2$. It follows that the vertical component of $g_t\gamma$ has length at least $\epsilon_\mathcal{K}/2$. In other words, $e^{-t} s \geq \epsilon_\mathcal{K}/2$. From this and \eqref{eq:B} we deduce the desired bound:
    \[
    k|w| \leq 1/s \leq \epsilon_{\mathcal{K}} e^{-t}/2.
    \]
    We conclude that $(J,R)$ has the desired properties if
    \[
        C = C(\mathcal{K}) \geq \max\{1,2,c_\mathcal{H},\epsilon_\mathcal{K}/2\}. \qedhere
    \]
\end{proof}
 
\subsection{Localization in the tautological plane}
\label{subsec:local}

In this section we prove that rigidity configurations can be applied to sharply constrain (candidate) eigenvalues along the 
unstable line in the tautological plane.

In the localization argument, it is necessary to estimate Hodge norms of the Poincar\'e duals of homology classes of
loops over compact sets of the moduli space. 

Given a translation surface $(X,\omega)$ and a closed piecewise geodesic curve $\gamma$ on it, denote by $[\gamma] \in H_1(X;\Z)$ its homology class, by $\#[\gamma] \in H^1(X;\Z)$ the corresponding Poincaré dual, and by $\ell_\omega(\gamma)$ the length of $\gamma$ with respect to the singular Euclidean metric induced by $\omega$ on $X$. The following estimate relating the Hodge norm of the Poincaré dual of a closed piecewise geodesic curve on a translation surface to its flat length will be useful for our purposes:

\begin{prop}
    \label{prop:hodge_bd}
    Let $\mathcal{H}^1$ be a stratum of unit area Abelian differentials. Then, for every compact set $\mathcal{K} \subseteq \mathcal{H}^1$, there exists a constant $C' = C'(\mathcal{K}) \geq 1$ such that for every translation surface $(X,\omega)\in \mathcal{K}$ and every closed piecewise geodesic curve $\gamma$ on $X$, the following estimate holds: 
    \[
    \|\#[\gamma]\|_X \leq C' \cdot  \ell_\omega(\gamma).
    \]
\end{prop}
\begin{proof} The statement can be derived from this result follows directly from \cite[Lemmas 3.2 and 5.5]{athfor}.
We sketch below a more direct argument for the convenience of the reader. The {\it stable norm} on $H_1(X, \R)$ induced by the flat metric of a translation surface $(X,\omega)$
is defined as 
$$
\Vert  h  \Vert_{\omega} = \inf \{ \sum_{i=1}^\sigma \vert r_i \vert \ell_{\omega}(\gamma_i) \,: \,  h= \Big[\sum_{i=1}^\sigma r_i \gamma_i \Big] \in H^1(X, \R)\}\,,
$$
where $\gamma_1, \dots, \gamma_{\sigma}$ are integral cycles.
Since $H_1(X, \R)$ is a finite dimensional vector space
the Poincar\'e dual of the Hodge norm  and the stable norm are
equivalent, that is, for all $(X, \omega)$,
$$
0< \inf_{h \in H_1(X, \R)}  \frac{\Vert h \Vert_ \omega} {\Vert \# h  \Vert_X} \leq \sup_{h \in H_1(X, \R)}  \frac{\Vert  h  \Vert_\omega} {\Vert  \# h  \Vert_X} < +\infty\,.
$$
In addition, both the Hodge norm and the stable norm depend continuously on $(X, \omega) \in \mathcal H^1$, hence for every
compact subset $\mathcal K \subset \mathcal H^1$, there exists 
a constant $C'(\mathcal{K})>0$ such that
$$
C'(\mathcal{K})^{-1} \leq \inf_{h \in H_1(X, \R)}  \frac{\Vert  h  \Vert_ \omega} {\Vert  \# h  \Vert_X} \leq \sup_{h \in H_1(X, \R)}  \frac{\Vert   h  \Vert_\omega} {\Vert \# h  \Vert_X} \leq C'(\mathcal{K})\,.
$$
In particular, for all $(X, \omega)\in \mathcal K$ and rectifiable loop $\gamma$, we have
$$
\Vert \# [\gamma] \Vert_X \leq  C'(\mathcal{K}) \Vert [\gamma]\Vert_\omega \leq  C'(\mathcal{K}) \ell_\omega(\gamma)\,.
$$
\end{proof}

Given a translation surface $(X,\omega)$ one can decompose its cohomology group 
\[
H^1(X;\R) = \Taut(X,\omega) \oplus \Bal(X,\omega),
\]
where $\mathrm{Taut}(X,\omega)$ is the tautological plane introduced in equation \eqref{eq:taut} and $\Bal(X,\omega)$ is the so-called `balanced' space, i.e., the symplectic complement (orthogonal) of the tautological plane with respect to the cup intersection form (Hodge inner product). Recall that we can further decompose
\[
\mathrm{Taut}(X,\omega) = \R \cdot \re (\omega) \oplus \R \cdot \im (\omega).
\]
In this context, denote by $\Bal_\mathcal{H}$ the bundle with fiber $\Bal(X, \omega)$ at each $(X, \omega) \in \mathcal{H}$,
and $\Taut^+_\mathcal{H}$ the unstable tautological subbundle, that is the 1 real dimensional subbundle whose fiber is 
$\R \im(\omega)$ at $(X, \omega)$.
\begin{gather*}
    \pi_\Bal \colon H^1_\R \to \Bal_\mathcal{H}, \\
    \pi_{\Taut^+} \colon H^1_\R \to \Taut^+_\mathcal{H},
\end{gather*}
the corresponding projections.

The following result will later be used in a parameter exclusion argument in the proof by contradiction of (2) implies (1) in Theorem \ref{thm:main1}.

\begin{prop} 
\label{prop:rig2_modified}
Let $\mathcal{H}^1$ be a stratum of unit area Abelian differentials over a moduli space of Riemann surfaces $\mathcal{M}$. Then, for every compact set $\mathcal{K} \subseteq \mathcal{H}^1$, there exist constants $C'' = C''(\mathcal{K}) \geq 1$ and $t_0 := t_0(\mathcal{K}) > 0$ with the following property. Let $(X,\omega) \in \mathcal{H}^1$ be any translation surface. Consider $t \geq t_0$ such that the following conditions hold:
\begin{enumerate}
\item $(X_t,\omega_t) \in \mathcal{K}$.
\item There exists $\vec{m} \in H^1(X_t; \Z)$ such that 
$d \left(g_t \cdot \alpha \im (\omega), \vec{m} \right) < 1/[2(C'')^2]$.
\end{enumerate}
Suppose $\alpha \in \R\setminus\{0\}$ has the property that   $d_\R(V\alpha,\Z)<1/C''$ for all $V > e^{t}/C''$ whenever $(X,\omega)$ has a $(V,H,1/C'',L)$-rigidity configuration.
Let $\beta_t \in \R$ be such that $\pi_{\Taut^+}(\vec{m}) = g_t \cdot \beta_t \im (\omega) \in H^1(X_t, \R)$. Then,
\[
d_\R(\alpha,\beta_t) \leq C'' \cdot \left\| \pi_\Bal(\vec{m}) \right\| \cdot e^{-t}.
\]
\end{prop}
\begin{proof}
For the rest of the proof fix a compact set $\mathcal{K} \subseteq \mathcal{H}^1$. Let $L_0 = L_0(\mathcal{K}) \geq 1$ and $C = C(\mathcal{K}) \geq 1$ be as in Proposition~\ref{prop:rigidity} and $C'= C'(\mathcal{K}) \geq 1$ be as in Proposition~\ref{prop:hodge_bd}. Define the constants $t_0 = \log L_0$ and
\[
C''=C''(\mathcal{K}) := \max\{20 C^2,C + 2C' +C'C, 20C\}.
\]
Now let $(X,\omega) \in \mathcal{H}^1$ be a translation surface whose vertical flow is ergodic and $\alpha,t$ be as in the statement of the proposition. 
In this setting we write
\[
\vec{m} = a \im (\omega_t) + b \im (\omega_t) + \pi_\Bal(\vec{m}).
\]
for some $a,b \in \R$. If $\beta:=\beta_t \in \R$ is such that $\pi_{\Taut^+}(\vec{m}) = g_t \cdot \beta \im (\omega)$, then
\begin{equation}
\label{eq:AA1}
a = \beta \cdot e^t.
\end{equation}
Condition (2) above ensures that
\begin{equation}
\label{eq:AA2}
|b| \leq 1/[2(C'')^2] \quad \text{and} \quad \|\pi_\Bal(\vec{m})\| \leq 1/[2(C'')^2].
\end{equation}
Consider now a $(V,H,\sigma,L)$-rigidity-configuration of $(X,\omega)$ for some parameters $V,H,\sigma, \allowbreak L > 0$, let $\gamma$ be a corresponding rigidity curve, and $g_t \gamma$ be its parallel transport to $X_t$. Denote by $[\gamma] \in H_1(X;\Z)$ the  homology class of $\gamma$ and by $g_t [\gamma]$ its parallel transport to $H_1(X_t;\Z)$. Letting $\cap$ be the pairing between cohomology and homology on $X_t$ we get
\[
\vec{m} \cap g_t[\gamma] = a \im (\omega_t) \cap g_t[\gamma] + b  \re (\omega_t) \cap g_t[\gamma] + \pi_\Bal(\vec{m}) \cap g_t [\gamma] .
\]
Notice that $\vec{m} \cap g_t[\gamma] \in \Z$. It follows that
\begin{equation}
\label{eq:A1}
d_\R(a \im (\omega_t) \cap g_t[\gamma],\Z) \leq |b| \cdot |\re (\omega_t) \cap g_t[\gamma]| + |\pi_\Bal(\vec{m}) \cap g_t [\gamma]|.
\end{equation}
Directly from the definitions we see that
\begin{equation}
\label{eq:A2}
\im (\omega_t) \cap g_t[\gamma] = e^{-t} \cdot V, \quad \re (\omega_t) \cap g_t[\gamma] = e^t \cdot H.
\end{equation}
By the Cauchy--Schwarz inequality
\begin{equation}
\label{eq:A3}
|\pi_\Bal(\vec{m}) \cap g_t [\gamma]| \leq \|\pi_\Bal(\vec{m})\| \cdot \| \#g_t[\gamma]\|
\end{equation}
By Proposition \ref{prop:hodge_bd},
\begin{equation}
\label{eq:A4}
\| \#g_t[\gamma]\| \leq C' \cdot \ell_{\omega_t}(g_t \gamma) = C' \cdot (e^{-t} \cdot  V + e^t \cdot H).
\end{equation}
Putting together \eqref{eq:AA1}, \eqref{eq:AA2}, \eqref{eq:A1}, \eqref{eq:A2}, \eqref{eq:A3}, and\eqref{eq:A4} we deduce
\begin{align}
\label{eq:C1}
d_\R(V \cdot \beta, \Z) &\leq  e^t \cdot H/[2(C'')^2] + C' \cdot \|\pi_\Bal(\vec{m})\| \cdot (e^{-t} \cdot V + e^t \cdot  H).
\end{align}
Consider now the subset $U = U(X,\omega,t,\vec{m}) \subseteq \R$ given by
\[
U := \{\rho \in \R \colon d(g_t \cdot\rho \im (\omega), \vec{m}) \leq 1/2C'' \}.
\]
By definition, $\alpha,\beta \in U$. Notice that if $\rho,\rho' \in U$, then the following holds:
\begin{align}
d_\R(\rho,\rho') &= d(\rho\im (\omega),\rho' \im (\omega)) = e^{-t} \cdot d(g_t\cdot\rho\im (\omega),g_t\cdot\rho' \im (\omega)) \nonumber  \\
&\leq e^{-t} \cdot d(g_t \cdot \rho\im (\omega),\vec{m}) + e^{-t} \cdot d (g_t \cdot\rho'\im (\omega),\vec{m}) \\
&\leq  \cdot e^{-t}/C''. \nonumber
\end{align}
In particular, we deduce the following bound on the diameter of $U \subseteq \R$:
\[
\mathrm{diam}_\R(U) \leq  e^{-t}/C''.
\]
Consider now a finite increasing sequence $L_1,\dots,L_k \in \R$ such that 
\begin{enumerate}
    \item $L_1 = e^t$.
    \item $L_k = e^t  / [2(C'')^2 \cdot \|\pi_\Bal(\vec{m})\| ]$.
    \item $L_i \geq 5C^2 \cdot  L_{i+1}/C'', \ \forall i =1,\dots,k-1$.
\end{enumerate}
By the choice of the constant $C''>0$ the bound $5 C^2 \cdot (1/C'') \leq 1/2$ holds, so finding such a sequence is always possible. As $t \geq t_0$ and $(X_t,\omega_t) \in \mathcal{K}$, Proposition~\ref{prop:rigidity} provides $(V_i,H_i,\sigma_i,L_i)$-rigidity-configurations $(J_i,R_i)$ of $(X,\omega)$ such that for every $i = 1,\dots,k$:
\begin{enumerate}
    \item $L_i/C \leq V_i \leq C \cdot L_i$.
    \item $\sigma_i \geq 1/C$.
    \item $H_i \leq C \cdot e^{-t}$.
\end{enumerate}
For every $i = 1,\dots,k$ denote by $\gamma_i$ a rigidity curve on $(X,\omega)$ corresponding to the rigidity configuration $(J_i,R_i)$ obtained above. Applying \eqref{eq:C1} to the rigidity curves $\gamma_i$ shows that for every $i=1,\dots,k$,
\begin{equation}
\label{eq:C2}
d_\R(V_i \cdot \beta, \Z) \leq  e^t \cdot H_i/[2(C'')^2] + C' \cdot \|\pi_\Bal(\vec{m})\| \cdot (e^{-t} \cdot V_i + e^t \cdot  H_i).
\end{equation}
Recall that $L_i \leq L_k = e^t  / [2(C'')^2 \cdot \|\pi_\Bal(\vec{m})\|]$ for every $i=1,\dots,k$. This together with the definition of $C''$, the bounds on the rigidity configurations $(J_i,R_i)$ presented above, and \eqref{eq:C2} shows that for every $i=1,\dots,k$,
\begin{equation}
\label{eq:D1}
d_\R(V_i \cdot \beta, \Z) \leq  1/(2C'').
\end{equation}
Notice now that, as $t \geq t_0$,
since  $(J_i, R_i)$ is, for every $i = 1,\dots,k$, a rigidity configuration with $V_i \geq e^t/C \geq e^{t_0}/C$ and
$\sigma_i \geq 1/C$, by the hypothesis on the eigenvalue $\alpha\not =0$ we derive that 
\begin{equation}
\label{eq:D2}
d_\R(V_i \cdot \alpha, \Z) \leq 1/(2C'').
\end{equation}
To summarize, we have proved that:
\begin{enumerate}
    \item $\alpha,\beta \in U \subseteq \R$.
    \item $\mathrm{diam}_\R(U) \leq  e^{-t} /C''$.
    \item $d_\R(V_i \cdot \alpha, \Z) \leq  1/(2C'')$ for every $i=1,\dots,k$.
    \item $d_\R(V_i \cdot \beta, \Z) \leq  1/(2C'')$ for every $i = 1,\dots,k$.
\end{enumerate}
We claim this implies that
\[
d_\R(\alpha,\beta) \leq 2 \cdot C \cdot C'' \cdot e^{-t} \cdot \|\pi_\Bal(\vec{m})\|.
\]
Indeed, proceeding by induction, one can show that, using the information up to step $i = 1,\dots,k$ constraints $\alpha$ and $\beta$ to belong to the same interval of length $1 /(C'' V_i)$ centered at some point in the lattice $\Z/{V_i}$; the bound $L_i \geq 5 C^2  L_{i+1}/C''$ guarantees that $V_{i+1} \leq (C''  V_i) /4$, making sure the intervals of \veigenvalues~are not getting split up at any step.  For $i=k$, since $ V_k \geq L_k/C$ and by the definition of $L_k$, we have that $\alpha$ and $\beta$ belong to the same interval of length
$$
\frac{1}{C''V_k} \leq \frac{C}{C''V_k} \leq \frac{C^2}{C''L_k}\leq 2 \cdot C\cdot C'' \cdot e^{-t} \cdot\|\pi_\Bal(\vec{m})\|\,,
$$
hence the claim follows immediately.
 \qedhere
\end{proof}

\subsection{Lattice repelling.} 
\label{subsec:repelling}

The following technical result follows directly from \cite[Proposition 5.1]{CF}. Roughly speaking, this result shows that if one gets fairly close to integral cohomology, one can push away from it after a bounded horocycle perturbation, thus breaking the Veech criterion and leading to weak mixing.

\smallskip
\noindent {\bf Notation}: We denote below, every $t,s \in \R$, $$(X_{t,s}, \omega_{t,s}) := g_t h_s (X,\omega) \in \mathcal{H}\,.$$

\begin{prop}
    \label{prop:CF_redux}
    There exist constants $0 < a < b < 1$ and $0 < \kappa < 1$ with the following property. Let $\mathcal{H}$ be a stratum of unit area Abelian differentials over a moduli space of Riemann surfaces $\mathcal{M}$, $(X,\omega) \in \mathcal{H}$, and $\mathcal{N} \subseteq \mathcal{H}$ be the $\mathrm{SL}(2,\R)$-orbit-closure of $(X,\omega)$. Then, there exist constants $t_0 = t_0(X,\omega) > 0$, $\sigma = \sigma(X,\omega) > 0$, and a compact set $\mathcal{K}_0 = \mathcal{K}_0(X,\omega) \subseteq \mathcal{H}$, such that for every $\mathrm{SL}(2,\R)$-invariant, strongly irreducible subbundle $\mathcal{F} \subseteq H_\R^1|_\mathcal{N}$ of the absolute cohomology bundle  over $\mathcal{N}$ with at least one positive Lyapunov exponent, there exists a constant $\tau = \tau(X,\omega,\mathcal{F}) > 0$ such that for every compact set $\mathcal{K} \supseteq \mathcal{K}_0$ there exists a constant $\gamma_0 = \gamma_0(X,\omega,\mathcal{F},\mathcal{K}) \in (0,1)$ with the following property.  Suppose that $s_0 \in [-1,1]$, that $t > t_0$, and that $\vec{v} \in \mathcal{F} \cap H^1(X, \R)$ are such that
    \begin{enumerate}
        \item $(X_{t, s_0}, \omega_{t, s_0}) \in \mathcal{K}$.
        \item $\gamma := \|g_t h_{s_0} \cdot \vec{v} \| \in (0,\gamma_0)$.
    \end{enumerate}
    Then there exists a set $\mathcal{S}_{t,s_0} = \mathcal{S}_{t,s_0}(X,\omega,\mathcal{F},\mathcal{K}) \subseteq [s_0 - e^{-2t}, s_0 + e^{-2t}]$ with
    \[
    \mathrm{Leb}(\mathcal{S}_{t,s_0}) > 2(1-\kappa^{-\sigma}\gamma^\sigma)e^{-2t}
    \]
    and such that for every $s \in \mathcal{S}_{t,s_0}$ there exists $\ell := \ell(X,\omega,\mathcal{F},\mathcal{K},s) > 0$ such that
    \begin{enumerate}
        \item $a |\log \gamma | \leq \ell \leq b |\log \gamma|$,
        \item $(X_{t+\ell,s}, \omega_{t+\ell,s}) \in \mathcal{K}$,
        \item $\|g_{t+\ell}h_s \cdot \vec{v}\| > e^{\tau \ell} \kappa \gamma$.
    \end{enumerate}
\end{prop}

Motivated by Lemma~\ref{lem:lattices} and Theorems~\ref{thm:galois_wright} and~\ref{thm:galois_filip}, we would like to apply Proposition \ref{prop:CF_redux} when $\mathcal{F}= \Gal(\mathcal{N})_\R$. This subbundle, although $\mathrm{SL}(2,\R)$-invariant, is not strongly irreducible. Nevertheless, following \cite[Theorem 1.2]{Fil16}, it can be decomposed as a direct sum of Hodge-orthogonal, $\mathrm{SL}(2,\R)$-invariant, strongly irreducible subbundles. Furthemore, Theorem \ref{thm:galois_filip} guarantees each component has at least one positive Lyapunov exponents. In particular, we deduce the following corollary.

\begin{cor}
    \label{cor:CF_redux}
    There exist constants $0 < a < b < 1$ and $0 < \kappa < 1$ with the following property. Let $\mathcal{H}^1$ be a stratum of unit area Abelian differentials over a moduli space of Riemann surfaces $\mathcal{M}$, $(X,\omega) \in \mathcal{H}$, and $\mathcal{N} \subseteq \mathcal{H}$ be the $\mathrm{SL}(2,\R)$-orbit-closure of $(X,\omega)$. Then, there exist constants 
    \begin{itemize}
    \item $t_0 = t_0(X,\omega) > 0$, \item $\sigma = \sigma(X,\omega) > 0$, \item $\tau = \tau(X,\omega) > 0$,
    \item $N =  N_\mathcal{N} \in \mathbb{N}$, and 
    \item a compact set $\mathcal{K}_0 = \mathcal{K}_0(X,\omega) \subseteq \mathcal{H}$, 
    \end{itemize} such that for every  compact set $\mathcal{K} \supseteq \mathcal{K}_0$, there exists a constant $\gamma_0 = \gamma_0(X,\omega,\mathcal{K}) \in (0,1)$ with the following property.  Suppose that $s_0 \in [-1,1]$, that $t > t_0$, and that $\vec{v} \in \Gal(\mathcal{N})_\R \cap \Bal(X,\R)$ are such that
    \begin{enumerate}
        \item $(X_{t, s_0}, \omega_{t, s_0}) \in \mathcal{K}$,
        \item $\gamma := \|g_t h_{s_0} \cdot \vec{v}\| \in (0,\gamma_0)$.
    \end{enumerate}
    Then there exists a set $\mathcal{S}_{t,s_0} = \mathcal{S}_{t,s_0}(X,\omega,\mathcal{K}) \subseteq [s_0 - e^{-2t}, s_0 + e^{-2t}]$ with
    \[
    \mathrm{Leb}(\mathcal{S}_{t,s_0}) > 2(1-\kappa^{-\sigma}{N^{\sigma/2}} \gamma^\sigma)e^{-2t}
    \]
    and such that for every $s \in \mathcal{S}_{t,s_0}$ there exists $\ell := \ell(X,\omega,\mathcal{K},s) > 0$ such that
    \begin{enumerate}
        \item $ a |\log \gamma|  \leq \ell \leq b |\log \gamma|$,
        \item $(X_{t+\ell,s},  \omega_{t+\ell,s}) \in \mathcal{K}$,
        \item $ \kappa^{-1} \gamma^{1-b}\geq \|g_{t+\ell}h_s \cdot \vec{v}\| > e^{\tau \ell} \kappa \gamma/N^{1/2} $.
    \end{enumerate}
\end{cor}
\begin{proof}
By Filip's semisimplicity theorem \cite{Fil16}, the bundle $\Gal (\mathcal N)_\R \cap \Bal(X,\R)$  over the invariant orbifold $\mathcal N$ is the direct, Hodge orthogonal, sum of strongly irreducible subbundles: there exists $N:=N_\mathcal{N} \in \N$
and strongly irreducible subbundles $\mathcal{F}_1, \dots, \mathcal{F}_N$ over $\mathcal{N}$ such that
$$
\Gal (\mathcal N)_\R \cap \Bal(X,\R) \vert \mathcal{N} =\mathcal{F}_1 \oplus \dots \mathcal{F}_N\,.
$$
Let $\{\pi_i \vert i\in \{1, \dots, N\}$ denote the projection
given by the above splitting. For any $\vec{v}\in \Gal (\mathcal N)_\R \cap \Bal(X,\R) $, there exists $j\in \{1, \dots, N\}$ such that 
$$
\Vert \pi_j (\vec{v}) \Vert = \max\{ \Vert \pi_i (\vec{v})\Vert : i\in \{1, \dots, N\}\}\,.
$$
Notice that as a consequence and by orthogonality 
$$
\gamma:=\Vert \vec{v}\Vert  \geq \gamma_j:= \Vert \pi_j (\vec{v}) \Vert \geq  \Vert \vec{v} \Vert /N^{1/2} =\gamma / N^{1/2}\,.
$$
We now apply Proposition~\ref{prop:CF_redux} to the vector
$\pi_j(\vec{v}) \in \mathcal{F}_j$, a strongly irreducible
subbundle with all Kontsevich--Zorich exponents strictly positive.

Suppose $\Vert \vec{v} \Vert < \gamma_0$, hence 
 $\Vert \pi_j(\vec{v}) \Vert =\gamma_j < \gamma_0$, 
so that by Proposition~\ref{prop:CF_redux} the following holds. For $t>t_0$ such that $(X_{t,s_0}, \omega_{t,s_0}) \in \mathcal{K}$ there exists a set $\mathcal{S}_{t,s_0} \subset [s_0-e^{-2t}, s_0+ e^{-2t}]$ such that 
\[
    \mathrm{Leb}(\mathcal{S}_{t,s_0}) > 2(1-\kappa^{-\sigma} \gamma_j^\sigma)e^{-2t} \geq 2(1-\kappa^{-\sigma} N^{\sigma/2}\gamma^\sigma)e^{-2t}
    \]
and such that for every $s \in \mathcal{S}_{t,s_0}$ there
exists $\ell \in \N$ such that $(X_{t+\ell,s}, \omega_{t+\ell,s})\in \mathcal{K}$, 
$$
a \vert \log \gamma_j \vert  \leq \ell \leq b \vert \log \gamma_j \vert \quad \text{ and } \quad \Vert g_{t+\ell} h_s \cdot \pi_j(\vec{v}) \Vert \geq e^{\tau \ell} \kappa \gamma_j\,.
$$
From the above bounds we derive that
$$
a \Big(\vert \log \gamma \vert -\frac{\log N}{2} \Big)   \leq \ell \leq b \vert \log \gamma_j \vert \quad \text{ and } \quad \Vert g_{t+\ell} h_s \cdot \vec{v} \Vert \geq e^{\tau \ell} \kappa \gamma/N^{1/2}\,.
$$
and\footnote{using the bound on the growth of the Hodge norm \cite[Section 2]{F02}}
$$
\Vert g_{t+\ell} h_s \cdot \vec{v} \Vert  \leq e^{\ell}  \Vert h_s\cdot {v} \Vert \leq  \kappa^{-1} \gamma^{-b} \Vert \vec{v}\Vert = \kappa^{-1} \gamma^{1-b}\,.
$$
By replacing the value of $\gamma_0$ and $a\in (0,1)$ from the Proposition \ref{prop:CF_redux}  by the values
$$
\gamma_0':=\min (\gamma_0, 1/ N)  \quad \text{ and }  a':= a/2\,,
$$
the statement holds and the proof is complete. 
\end{proof}

\subsection{Proof (2) implies (1) in Theorem \ref{thm:main1}}

Throughout this section,  the compact set $\mathcal{K}' \subset \mathcal{H}^1$ is chosen to be large enough so that all of the above results hold for it (in particular $\mathcal{K}\supset \mathcal{K}_0$ where $\mathcal{K}_0$ is as in Corollary \ref{cor:CF_redux}) and we choose $\mathcal{K} $ a larger compact set so that $h_s(X,\omega) \in \mathcal{K}$ whenever $|s|\leq 1$ and $(X,\omega) \in \mathcal{K}'$. The parameter $\hat{\epsilon}>0$ will be chosen small enough that a number of smallness conditions in the proof hold. As all of these conditions hold for all $\epsilon>0$ small enough, this does not cause issues. Given these terms the Teichm\"uller geodesic flow time parameter is chosen large enough (see \eqref{eq:key set}).

\begin{lem} 
\label{lem:set A}
Let $\alpha: \mathcal{S} \to \R\setminus \{0\}$ be a measurable eigenvalue section (as in Proposition~\ref{prop:selection}).
There exists a compact subset  $K_{\mathcal{S}} \subset \mathcal{S}$ with $|K_{\mathcal{S}}\setminus \mathcal{S}|<\frac \epsilon 9$ such  that $\alpha|_{K_{\mathcal{S}} }$ is continuous. For any compact set $\mathcal{K} \subset \mathcal{H}^1$ and for any $\epsilon>0$, there exists  $t_*:=t_*(\epsilon, \mathcal{K})\geq 2$ such that 
\begin{equation}\label{eq:key set} A(\epsilon,\mathcal{K}, t_*):=\{s \in K_\mathcal{S}: \alpha(s) \text{ is a }(\epsilon,\mathcal{K}, t_*)-\text{\veigenvalue}\}
\end{equation} 
has positive measure. 
\end{lem}
\begin{proof}
The existence of the compact set $K_\mathcal{S}$ of large
measure and such that $\alpha\vert K_\mathcal{S}$ is continuous
follows from Proposition~\ref{prop:selection} by Lusin's theorem. 
By Lemma \ref{lem:candidate mble}, for every $n\in \N$, the set
$$
A_n := \{s \in K_\mathcal{S}: \alpha(s) \text{ is a }(\epsilon,\mathcal{K}, n)-\text{\veigenvalue}\}
$$
is measurable. By Lemma \ref{lem:candidate pos}, since $A_n \subset A_{n+1}$ for all $n\in \N$, and 
$$
K_\mathcal{S} = \cup_{n\in \N\setminus\{0,1\}} \,A_n\,,
$$
by measure continuity there exists $n_*\in \N\setminus\{0,1\}$ such that $A_n$ has positive measure. We then let 
$t_*=n_*$ and $A= A_{n_*}$ and the argument is concluded.
\end{proof}

\begin{prop}\label{prop:end intermediate}
    Let us assume that there exists a positive measure set of $\theta\in [0, 2\pi)$ so that the flow in direction $\theta$ on $(X,\omega)$,  is not weakly mixing. Let $\mathcal{K},\mathcal{K}'
    $ be the compact sets defined at the start of this section. There exists $\epsilon_0(\mathcal{K}')>0$ such that if $0<\epsilon <\epsilon_0(\mathcal{K}')$ there exists $t_0>0$,
    $s_0\in [-1,1]$ with the property that 
    \begin{enumerate}
    \item  $(X_{t_0,s_0},\omega_{t_0,s_0}) \in \mathcal{K}'$,
    \item the vertical flow $\Phi_{s_0}$ on $(X_{0,s_0}, \omega_{0,s_0})$ has an eigenvalue $\alpha_0$,
    \item\label{conc:large} the set $$
    \begin{aligned}\mathcal E:= \{s\in& [s_0-e^{-2{t_0}},s_0+e^{-2t_0}] : \text{the vertical flow $\Phi_{s}$ on $(X_{0,s}, \omega_{0,s})$} \\ &\text{has a $(\epsilon,\mathcal{K}, t_0)$-\veigenvalue~$\alpha_s$ in $B(\alpha_0,\epsilon e^{-t_0})$} \}
    \end{aligned}
    $$
    has Lebesgue measure at least $2(1-\epsilon)e^{-2t_0}$,
    \item there is $\vec{n} \in  \Lambda(\mathcal{N})\setminus \{\vec{0}\} $ such that, for all $s\in \mathcal{E}$,
    $$ d \Big(g_{t_0} \cdot \alpha_s \im(h_{s}\omega), g_{t_0} h_s \cdot \vec{n}\Big)<{\epsilon}\,.$$

    \end{enumerate}
\end{prop}
\begin{proof}
Let $\alpha: \mathcal{S} \to \R\setminus \{0\}$ be a measurable eigenvalue section (as in Proposition~\ref{prop:selection}) and let  $K_\mathcal{S} \subset \mathcal{S}$ and 
$A:= A(\hat{\epsilon},\mathcal{K}, t_*)=\{ s \in K_\mathcal{S}: \alpha(s) \text{ is a }(\hat{\epsilon},\mathcal{K}, t_*)-\text{\veigenvalue} \}$~be the sets of Lemma \ref{lem:set A}. The parameter $\hat{\epsilon}>0$ is to be chosen later.
We choose $s_0\in (-1,1)$ to be a density point of $A$ so that $(X_{\ell,s_0},\omega_{\ell,s_0})\in \mathcal{K}$ for a positive density set of $\ell$, which is true for almost every $s$ by \cite{CE15}, and define $\alpha_0:=\alpha(s_0)$. 

We now choose $\delta>0$ so that 
\begin{itemize}
\item $(s_0-\delta,s_0+\delta)\subset [-1,1]$,
\item for all $0<r\leq \delta$ we have 
$$|(s_0-r,s_0+r)\cap A|>2r(1-\frac \epsilon 9),$$
\item  for all $0<r\leq \delta$, $$\text{\rm diam}\Big(\alpha(K_\mathcal{S} \cap (s_0-r,s_0+r)\Big)<\frac 1 9 e^{-t_*}\,.$$
\end{itemize}
For any $t_0>\max\{t_*,-\log(\delta)\}$  chosen so that $(X_{t_0,s_0}, \omega_{t_0,s_0}) \in \mathcal{K}'$ we have

\medskip
\noindent
\textbf{{Claim:}} Let $C=C(\mathcal{K})>0$ denote the constant of Proposition~\ref{prop:rigidity}.  If $\hat{\epsilon}<\frac 1 {9C^2}$,  $s \in [s_0-e^{-2t_0},s_0+e^{-2t_0}]$, $\alpha$ is an $(\hat{\epsilon},\mathcal{K},t_*)$-\veigenvalue~of the vertical flow of $h_{s}(X,\omega)$ and $|\alpha-\alpha_0|<\hat{\epsilon} e^{-t_*}$, then for every $\sigma \geq 1/C$,
\begin{equation}\label{eq:stuck close}|\alpha-\alpha_0|<\frac{2C^2}{\sigma} \hat{\epsilon} \, e^{-t_0}. 
\end{equation}

\begin{proof}[Proof of Claim]
Observe that if $(J,R)$ is a $(V,H,\sigma,L)$-rigidity configuration for $(X',\omega')$, then we can construct a $(V,H',\sigma-\vert\ell\vert V^2,L)$-rigidity configuration with $H'\in [H-\vert\ell\vert V,H+\vert\ell\vert V]$ of $h_{\ell}(X',\omega')$ for all $|\ell| < \sigma/(2V^2)$.\footnote{Indeed, consider  $h_\ell R$ and $h_\ell R^*$, the images of $R,\, R^*$ under the affine homeomorphism given by $h_\ell$. These are are parallelograms with horizontal sides and the other sides have holonomies  $(\ell V,V)$. 
 Each horizontal segment of $R^*$ has length $\frac{\sigma}{V}$ and so if $|\ell|<\frac{\sigma}{V^2}$ then a portion of each horizontal segment of length $\frac{\sigma}{V}-|\ell|V$ flows under the vertical flow into the top segment of $h_\ell R^* \subset h_\ell J$. Thus we obtain a rigidity configuration on $h_\ell(X',\omega')$ where the relevant $R^*$ has measure  
 $\sigma -|\ell| V^2$ 
 }

As a consequence, if $(J,R)$ is a $(V,H,\sigma,L)$-rigidity configuration for $h_{s_0}(X,\omega)$, since $\vert s-s_0\vert < e^{-2t_0}$, under the assumption that $V<e^{t_0}\frac{\sigma}2$ we have a rigidity configuration for $h_{s}(X,\omega)$. Considering a sequence of such $V_1<V_2<....<V_k$ with consecutive ratio at most $C^2$, $V_1\in [e^{t_*},C^2e^{t_*}]$ and $V_k\in [\sigma C^{-2}e^{t_0},\sigma e^{t_0}]$. Now because $\hat{\epsilon}<\frac{1}{9C^2}$, if $$d_\R(\alpha' V_i,n_i) \,,  d_\R(\alpha'' V_i,n_i)<\hat{\epsilon}$$  and $$d_\R(\alpha' V_{i+1},\Z)\,, d_\R(\alpha'' V_{i+1},\Z)< \hat{\epsilon}$$ then there exists $n_{i+1}$ so that $$d_\R(\alpha' V_i,n_{i+1}), d_\R(\alpha''V_i,n_{i+1})<\hat{\epsilon}\,.$$
From our assumptions, by finite iteration of the above procedure, we derive that 
\begin{equation}\label{eq:pinned}|\alpha -\alpha_0 |\leq  
2(V_k)^{-1}\epsilon\leq  2\frac{C^2}{e^{t_0}\sigma}\hat{\epsilon}
\end{equation}
\end{proof}
From the claim (with $\sigma=1/C$) we obtain the first three conclusions of the proposition, if $\hat\epsilon$ is small enough to apply the claim and such that 
$$
 \frac{2 C^2}{\sigma} \hat \epsilon \leq 2 C^3 \hat \epsilon < \epsilon\,. 
$$

For $\alpha_0$ an eigenvalue of the vertical flow on $(X_{0,s_0}, \omega_{0,s_0})$ and by Lemma~\ref{lem:lattices},
since $\alpha_{s_0} \im(h_{s_0}\omega) \in \Gal(\mathcal{N})_\R$, for $t_0>0$ sufficiently large, under the assumption that $(X_{t_0,s_0}, \omega_{t_0,s_0})\in \mathcal{K}'$ there is $\vec{n}\in \Lambda(\mathcal{N})$ such that  
$$d \Big(g_{t_0} \cdot \alpha_{s_0} \im(h_{s_0}\omega), g_{t_0} h_{s_0} \cdot \vec{n}\Big)<{\hat{\epsilon}}.$$
Since $\im(h_s\omega)= \im(\omega)$ for all $s\in \R$, we have that, by the above Claim,
$$
d \Big(g_{t_0} \cdot \alpha_{s_0} \im(h_{s_0}\omega), 
g_{t_0} \cdot \alpha_{s} \im(h_{s}\omega\Big) \leq 
 e^{t_0} d_\R (\alpha_s, \alpha_0) \leq 2C^3\hat{\epsilon}\,.
$$
We have the statement of (4)  if we choose $\hat{\epsilon}= \epsilon/(1+2C^3)$, hence the argument is complete.
 \end{proof}

We now prove that (2) implies (1) by contradiction, assuming that (2) is true and not (1) is true.

\begin{proof} 
We obtain $s_0$ as in the proof of 
Proposition \ref{prop:end intermediate} for $\epsilon>0$ to be chosen based on smallness conditions to come. 
By  
Proposition~\ref{prop:end intermediate}, item  (4), 
there is a single integer vector $\vec{n} \in  \Gal(\mathcal{N})_\R\cap H^1 (X, \Z)$ so that $$ g_{t_0} h_s \cdot \Big(\vec{n}-  \alpha_s \im(\omega) \Big) \in H^1(X_{t_0,s}, \R)$$ is small for all $s \in \mathcal {E} \subset [s_0-e^{-2t_0},s_0+e^{-2t_0}]$. Because $\vec{n}$ is a non-zero integer vector and by hypothesis $\Taut(X,\omega)$ contains no integer points, but it is contained in $\Gal(\mathcal{N})_\R$ which is defined over $\Q$, it follows that $$ \pi_{\Bal}(\vec{n})=  \pi_{\Gal(\mathcal{N})} \pi_{\Bal}(\vec{n}):=\vec{v}\neq 0\,.$$
By Filip's Theorem~\ref{thm:galois_filip}, $\vec{v}$ is contained in the sum of strongly irreducible subbundles of $\Bal \cap \Gal(\mathcal{N})$, each with a positive top Lyapunov exponent.

Because $0<\| g_{t_0}h_{s_0}\cdot\vec{v}\|< \epsilon\ll 1$, by Proposition \ref{prop:end intermediate}, item (4), we can now iteratively apply Corollary~\ref{cor:CF_redux} and Proposition~\ref{prop:rig2_modified}.

Let $u_0= v_0=\Vert g_{t_0}h_{s_0}\cdot \vec{v}\Vert <\uzero$ and 
$w_{0}=2e^{-2t_0}$. Recursively let 
\begin{equation}
\label{eq:recursion}
\begin{aligned} 
u_{i+1}&= \frac{\kappa}{N^{1/2}} \cdot u_i^{1-a\tau}\,, \\
v_{i+1}&= \kappa^{-1} u_{i+1}^{1-b}\,,\\
w_{i+1}&=(1- 2\kappa^{-\sigma} N^{-\sigma/2} u_{i+1}^{\sigma}) w_i \,,
\end{aligned}
\end{equation}
for $i=0,\dots, m-1$, where  
$$
m=\min\{i:u_i\geq \um \}\,.
$$
We inductively prove the existence of families of intervals $\mathfrak{I}_i=\{J_{i,j}\}_{j=1}^{n_i}$ so that 
\begin{enumerate}[label=(\alph*)]
\item $|\bigcup_{J\in \mathfrak{I}_i}J| \geq w_{i}$, 
\item For each $J_{i,j}$ there exists $s_{i,j}$ and $r_{i,j}$ so that $$J_{i,j}=[s_{i,j}-e^{-2(t_0+r_{i,j})},s_{i,j}+e^{-2(t_0+r_{i,j})}]\,,$$
\item $(X_{i,j}, \omega_{i,j}):=g_{t_0+r_{i,j}}h_{s_{i,j}}(X,\omega) \in \mathcal{K}{\color{brown}'}$,
\item $v_{i}>\Vert g_{t_0+r_{i,j}}h_{s_{i,j}}\cdot \vec{v}\Vert  \geq u_i$,
\item Each $s \in [s_0-e^{-2t_0},s_0+e^{-2t_0}]$ is in at most two intervals $J \in \mathfrak{I}_i$. 
\end{enumerate}

The base case is $\mathfrak{I}_0=\{[s_0-e^{-2t_0},s_0+e^{-2t_0}]\}$. We now do the inductive step. We will first apply Corollary~\ref{cor:CF_redux} to build a set of intervals $\mathfrak{I}_{i+1}'$, satisfying $(a)-(d)$. With this in hand, we can throw out some of the intervals to obtain $\mathcal{I}_{i+1}$ which satisfies (e) as well and so that 
$\bigcup_{J \in \mathfrak{I}_{i+1}'}J=\bigcup_{J \in \mathfrak{I}_{i+1}}J$. Indeed any set of intervals can be refined to a set of intervals with the same union so that each point belongs to at most two of the intervals. 

\smallskip
\noindent \textbf{Constructing $\mathfrak{I}_{i+1}'$}: Let $J=[s'-e^{-2(t_0+r')},s'+e^{-2(t_0+r')}] \in \mathfrak{I}_j$ be given. If 
\begin{equation}\label{eq:lazy} \Vert g_{t_0+r_{i,j}}h_{s_{i,j}}\vec{v}\Vert \geq u_{i+1} \,, \end{equation}
we put $J \in \mathfrak{I}'_{i+1}$ (with the same $s$ and $r$ parameters). 

Otherwise we apply Corollary \ref{cor:CF_redux} with $s_0=s'$, $t_0=t_0+r'$ and $$\gamma=\Vert g_{t_0+r_{i,j}}h_{s_{i,j}} \cdot\vec{v}\Vert<u_{i+1}\,.$$ We take the set of intervals given by the proposition and add them to $\mathfrak{I}_{i+1}'$.

Now \eqref{eq:lazy} and Corollary \ref{cor:CF_redux} first and third conclusion give conclusions (b)-(d). Conclusion (d) follows from the fact that the previous $$\Vert g_{t_0+r_{i,j}}h_{s_{i,j}}\cdot\vec{v}\Vert  \in [u_i, u_{i+1}]\,,$$ and conclusion 2. of Corollary \ref{cor:CF_redux}. 
  
The conclusion on the size of the measure follows inductively from considering the complement of the $\mathcal{S}_{i,j}:=\mathcal{S}_{t,s}$
with $t= t_0+ r_{i,j}$ and $s=s_{i,j}$. The largest these complements can be is if $\gamma$ is maximal which is when it is equal to $u_{i+1}$ and, by property (e), which is satisfied by the family of intervals $\mathcal I_i$, we obtain the coefficient $2$ in the definition of $w_{i+1}$. 

We now make conclusions based on this: 
\begin{lem}For all $\eta_0>0$ there exists $\epsilon_0>0$ so that if $\epsilon<\epsilon_0$ then $$w_m>(1-\eta_0)w_0\,.$$
\end{lem}
\begin{proof} 
From the recursion one can check that the $u_j$ grow exponentially for $\epsilon$ small enough. The desired bound follows from the definitions and this observation. More precisely, one can check that for every $D > 1$ there exists $\epsilon_D > 0$ such that for $0 <\epsilon < \epsilon_D$, the following bound holds:
\[
w_m \geq \prod_{j=-\infty}^{-1}(1-D^j)w_0. \qedhere
\]
\end{proof}

\begin{lem}\label{lem:eval connect_bis} For any $\eta_1>0$ there exists $\epsilon_1>0$ so that whenever $\uzero <\epsilon_1$ then
$$d( g_{t_0+r_m} \cdot \alpha_s \im(h_s\omega), g_{t_0+r_m} h_s \cdot \vec{n})<\eta_1 \,,$$
for all $s \in J \cap \mathcal{E}$ where $J=[s'-e^{-2(t_0+r_m)},s'+e^{-2(t_0+r_m)}]\in \mathfrak{I}_m$. 
\end{lem}

\begin{proof}We consider the sequence of strictly decreasing
intervals $I_1,I_2,...,,I_{m'}$ such that 
\begin{itemize}
\item $s \in I_\ell$ for all $\ell=1, \dots, m'$, \,;
\item there exists a strictly increasing  function 
$\psi:\{1,\dots,m'\} \to \{0,\dots,m\}$ 
with the property that  $\psi(1)=0$ and $\psi(m')=m$
such that
$$I_\ell:=[s_\ell-e^{-2(t_0+r_\ell)},s_\ell+e^{-2(t_0+r_\ell)}] \in \mathfrak{I}_{\psi(\ell)}\,, \quad \text{for all } \ell=1, \dots, m'\,;$$
\item {$I_{\ell+1}$ is obtained by applying Corollary \ref{cor:CF_redux} to $I_{\ell}$}
(and in particular $I_\ell \neq I_{\ell+1}$) for $\ell=1, \dots, m'-1$. 

 This gives that 
 \begin{equation}\label{eq:flow short} 0< r_{\ell+1}-r_{\ell}\leq 
 -b\log \Vert g_{t_0+r_{\ell}}h_{s_{\ell}} \cdot \vec{v}\Vert  \,.\end{equation}
\end{itemize} 

Let $\beta_{s,\ell}$ be the unique real number such that
$$
\pi_{\Taut^+}( g_{t_0+r_\ell} h_{s}\cdot \vec{n}) =  g_{t_0+r_l} \cdot \beta_{s,\ell} \im(h_s\omega)\,.
$$
Since $\alpha_s$ is a \veigenvalue,~Proposition \ref{prop:rig2_modified} implies that  for $t_0$ sufficiently large,
under the assumption that $d(g_{t_0+r_l}\cdot \alpha_s \im(h_s\omega), g_{t_0+r_\ell} h_s \cdot \vec{n}) < 1/[2(C'')^2]$,
$$
d_\R (\alpha_s, \beta_{s, \ell}) \leq C'' \Vert \pi_\Bal ( g_{t_0+r_\ell} h_s \cdot \vec{n}) \Vert e^{-(t_0+ r_\ell)} =
C'' \Vert  g_{t_0+r_\ell} h_s \cdot \vec{v} \Vert e^{-(t_0+ r_\ell)} \,.
$$
By taking into account that the tautological subbundle is the direct (orthogonal) sum of the unstable line subbundle $\Taut^+$ and of the stable line subbundle, and by the bound in formula \eqref{eq:flow short} we derive that there exists a constant $C'''>0$  (one can take $C'''= e C'' +1$) such that
$$
d(g_{t_0+r_{\ell+1}} \cdot \alpha_s \im(h_s\omega), g_{t_0+r_{\ell+1}} h_s \cdot \vec{n}) \leq C'''  \Vert  g_{t_0+r_\ell} h_s \cdot \vec{v} \Vert^{1-b} 
+ \uzero e^{-r_{\ell+1}}\,.
$$
Since by construction (and by the ``trivial'' upper bound on the growth of the Hodge norm) we have
$$
\label{eq:dist_est2}
\begin{aligned}\Vert  g_{t_0+r_\ell} h_s \cdot \vec{v} \Vert^{1-b} \leq  v_{\psi(\ell)}^{1-b}  &= (\kappa^{-1} u_{\psi(\ell)}^{1-b} )^{1-b}
\\ &\leq  \kappa^{-(1-b)} \Big((\um)     ( \um)^{-b}  \Big )^{(1-b)^2} = (\frac{1}{\kappa})^{1-b} (\um )^{(1-b)^3} \,,
\end{aligned}
$$
by choosing $\epsilon >0$ so small that 
$$
C'''(\frac{1}{\kappa})^{1-b} (\um)^{(1-b)^3} + \uzero <  1/[2 (C'')^2]\,,
$$
the required assumptions are verified by induction for all $\ell =1, \dots, m'$  starting with
$$
d(g_{t_0}\cdot \alpha_s \im(h_s\omega), g_{t_0} h_s \cdot \vec{n}) < \epsilon <1/[2(C'')^2]\,,
$$
a condition that  holds by Proposition \ref{prop:end intermediate} for all $s\in \mathcal{E} \subset [s_0-e^{-2t_0}, s_0+e^{-2t_0}]$, hence the desired conclusion follows as the last step of the finite induction.
\end{proof}

We now complete the proof.
First observe that there exists $ s_* \in \bigcup_{J \in \mathfrak{I}_m}J \cap \mathcal{E}$. Indeed, if $\uzero<1/2$ is small enough, then $w_{m}>\frac 1 2 e^{-2t_0}$ and since by Proposition \ref{prop:end intermediate} the Lebesgue measure of the $\mathcal{E}$ is at least $2(1-\epsilon)e^{-2t_0}$, we have such a $s_*$. 
Let $J=[s'-e^{-2(t_0+r')},s'+e^{-2(t_0+r')}]\in \mathcal{I}_m$ be the interval $s_*$ is in. 

By Proposition~\ref{prop:end intermediate}, item (4),   
there is $\vec{n} \in  \Lambda(\mathcal{N})\setminus \{\vec{0}\} $ such that, for all $s\in \mathcal{E}$,
\begin{equation} 
\label{eq:upper}
    d \Big(g_{t_0} \cdot \alpha_s \im(h_{s}\omega), g_{t_0} h_s \cdot \vec{n}\Big)<{\epsilon}\,.
\end{equation}    
On the other hand, we have that, for $s=s_*$,
\begin{equation}\label{eq:lower}
10 \uzero <\inf \Big\{ d\Big(g_{t_0+r'} h_s \cdot \alpha \im (\omega), g_{t_0+r'}h_s \cdot \vec{n}\Big):\alpha \in \R \Big\}.
\end{equation}


To see this, by our procedure and choice of $m$,  we have
$$\Vert g_{t_0+r'}h_{s'} \cdot \vec{v}\Vert \geq u_m \geq 
\um\,,$$
hence for $s=s_*$ (a horocycle parameter value possibly different from $s'$)
$$
\begin{aligned}
d\Big(g_{t_0+r'} h_ \cdot \alpha \im (\omega), g_{t_0+r'}h_s \cdot \vec{n}\Big) &\geq \kappa d\Big(g_{t_0+r'} h_{s'} \cdot \alpha \im (\omega), g_{t_0+r'}h_{s'} \cdot \vec{n}\Big) \\ & \geq \kappa \Vert g_{t_0+r'}h_{s'} \cdot \vec{v}\Vert \geq \kappa \cdot\um = 10 \uzero\,.
\end{aligned}
$$
 
Equations \eqref{eq:upper} and \eqref{eq:lower} are in contradiction, which completes the proof.
\end{proof}

\section{Weak mixing for rational polygons}
In this section we characterize rational polygons whose billiard flow is weakly mixing on Lebesgue almost every invariant surface.
We introduce the following
\begin{defin} A translation surface or a rational polygonal billiard has the \emph{weak mixing property} if its directional
flow is weakly mixing for Lebesgue almost all directions.
\end{defin}

Unfoldings of  polygons give cyclic covers of the Riemann sphere. We describe 
below the relevant part of the construction following M. Mirzakhani and A. Wright \cite[\S 6]{MW18}. 

\medskip
Let $\theta_1, \dots, \theta_n \in (0, \pi) \cup (\pi, 2\pi)$ such that $\theta_i/\pi \in \Q$ for all $i\in \{1, \dots, n\}$
and $\sum_{i=1}^n \theta_i = (n-2)\pi$. Let $k$  be the least common denominator of the finite set $\{\theta_i/\pi \vert  i=1, \dots, n\}$
and set $q_i= k \theta_i/ \pi$. For each tuple $z = (z_1, \dots, z_n) \in \C^n$ of distinct complex numbers, let
$X_z$ to be the normalization of the plane algebraic curve
$$
y^k = \prod_{i=1}^n (z-z_i)^{q_i} \,.
$$
This Riemann surface has an automorphism $T$ given by 
$$T(y, z) = (\xi y, z),  \quad \text{where }\xi= \exp(2\pi i/k)\,.$$
The automorphism $T$ generates the deck group for the covering map $(y,z) \to z$. The surface $X_z$ is a cyclic
cover of ${\mathbb P}^1(\C)$.

The differential 
$$
\omega_z = \prod_{i=1}^n (z-z_i)^{\frac{\theta_i}{\pi} -1} dz  =  y  \prod_{i=1}^n (z-z_i)^{-1} dz 
$$
is holomorphic on $X_z$ and such that (see \cite[Lemma 6.1]{MW18})
$$
T^\ast (\omega_z) = \xi \cdot \omega_z = \exp(2\pi i/k) \omega_z \,.
$$
 
\begin{lem} \cite[ Lemma 6.1]{MW18}, \cite{DT02} For any polygon with interior angles $\theta_1, \dots, \theta_n \in (0, \pi) \cup (\pi, 2\pi)$  as above, there exists $z=(z_1, \dots, z_n) \in \C^n$ such that its unfolding is the cyclic cover  $(X_z, \omega_z)$.
\end{lem}

\begin{lem} A rational polygonal billiard has the weak mixing property if and only if its unfolding has the weak mixing property.
\end{lem} 
\begin{proof}
The directional flows of a rational polygonal billiards are a subset of the directional flows of its unfolding, hence a rational
polygonal billiard has the weak mixing property if its unfolding has it. 

Conversely,  since directional flows of a rational polygonal billiards form a subset of positive measure (a non-trivial closed interval)
in the set of directional flows of its unfolding, it follows from the main Theorem \ref{thm:main1} that the unfolding has the weak mixing property if the rational polygonal billiard has it.
\end{proof}

 Let $P$ be a polygon with angles $\theta_1, \dots, \theta_n \in (0, \pi) \cup (\pi, 2\pi)$ such that $\theta_i/\pi \in \Q$ for all $i\in \{1, \dots, n\}$ and $\sum_{i=1}^n \theta_i = (n-2)\pi$. Let $k$  be the least common denominator of the finite set $\{\theta_i/\pi \vert  i=1, \dots, n\}$.
 
\begin{cor} \label{cor:gen_case}  If $k \not =1,  2, 3, 4, 6$, then  $P$ has the weak mixing property.
\end{cor} 
\begin{proof} If $P$ does not have the weak mixing property, neither does its unfolding $(X_z, \omega_z)$, hence by the main Theorem \ref{thm:main1}, the unfolding has an integer point (that is, an element of $H^1(X_z, \Z)$) in its tautological plane 
$H^1_{taut}(X_z, \R)=  \R\cdot \text{\rm Re}(\omega_z) +  \R\cdot\im (\omega_z)$.
Since the automorphism $T$ leaves $\omega_z$ invariant, it also leaves the tautological plane invariant and acts in it by a
rotation of angle  $2\pi/k$. Since $T: H^1(X_z, \Z) \to H^1(X_z,\Z)$, that is, it preserves the integer lattice, it follows
that the tautological plane contains a lattice with a $k$-fold rotational symmetry, with $k \not =1,  2, 3, 4, 6$, which is a contradiction since it is well-known that no such lattice exists (as it can be proved as an exercise).
\end{proof} 

It remains to examine the cases of $k=1, 2, 3, 4, 6$. 

\medskip
The case $k=1$ cannot occur by the condition $\sum_{i=1}^n \theta_i =(n-2) \pi$, since not all angle $\theta_i$ can be integer
multiples of $\pi$. 

\medskip
The cases $k=2, 3, 4, 6$ include all {\it almost integrable} polygons according to E. Gutkin's definition \cite{Gu84} and
all polygons of $\Delta$-class according to the definition of  Gutkin and A. Katok \cite{GK88} (see \cite[Defs. 2 and 3]{GK88}). 

\begin{defin}  \label{def:a_i} \cite{Gu84}, \cite{Gu86}, \cite{GK88} A polygon $P$ is called \emph{almost integrable} if the group $G_P$ generated by reflections in the sides of $P$ is a discrete subgroup of the group $E(\R^2)$ of rigid motions of the plane or, equivalently, if it is homothetical to a polygon drawn on the lattice generated by the unfolding of a completely integrable polygon.  Given a completely integrable polygon $\Delta$, a polygon is called of \emph{$\Delta$-class} if its sides are parallel to the lines of the lattice generated by the unfolding of $\Delta$.
\end{defin}
\begin{rem}  By the expression ``polygon drawn on the lattice generated by the unfolding of a completely integrable polygon'' it is meant a polygon whose boundary is a subset 
of the graph (grid) in $\R^2$ obtained by applying the group of reflections of an integrable polygon to the polygon itself.
\end{rem}

Gutkin and A. Katok proved in \cite{GK88} that, for any set of  polygons of $\Delta$ class with {\it fixed combinatorics} 
 (that is, fixed number of connected components of the boundary, fixed number of vertices and angles at each vertex) and sufficiently many edges (at least $5$ when $\Delta$ is a rectangle and at least $4$ when $\Delta$ is a triangle), for any irrational angle $\theta$ there is a $G_\delta$ dense set of polygons with weakly mixing directional flow in direction $\theta$. We will completely characterize below almost integrable polygons with the weak mixing property. In fact, we prove a strong version of the statement ``reasonably expected'' by Gutkin
 (see \cite{Gu84}, page 570) that almost all polygons which are not almost integrable have the weak mixing property.

\medskip
The case $k=2$ corresponds to polygons with vertical/horizontal edges, and we treat this later.

\begin{cor}  \label{cor:bct} If $k =3, 4, 6$, a polygon $P$ does not have the weak mixing property, if and only if its unfolding $(X_z, \omega_z)$  is a branched cover of a flat torus. This torus is rectangular if $k=4$, and hexagonal if $k=3$ or $k=6$.
\end{cor} 
\begin{proof}  In this case the automorphism $T$ acts on the tautological plane as a rotation of order $3$, $4$ or $6$. By Theorem
\ref{thm:main1} if $P$ does not have the weak mixing property, then the tautological plane contains a non-zero integer vector. Since
$T$ preserves the integer lattice $H^1(X_z,\Z)$ and has no non-trivial real eigenvectors, it follows that the tautological plane contains two linearly independent
integer vectors, hence it contains a sublattice of the integer lattice.  This sublattice is rectangular for $k=4$ and hexagonal if 
$k=3$ or $k=6$. In other terms, there exists a matrix $A \in GL(2, \R)$
such that  
$$
A \begin{pmatrix} \text{\rm Re}  (\omega_z) \\\im  (\omega_z) \end{pmatrix} \in H^1(X_z, \Z \oplus i \Z)\,.
 $$
Let then $\Hol: X_z \to \C$ be the multi-valued map defined as (for a given regular point $x_0 \in X_z$) 
$$
\Hol_z(x)  = \int_{x_0}^x   \omega_z \,, \quad \text{for all } x \in X_z\,.
$$
The map  $\Hol_z$ is well-defined modulo the lattice $A^{-1} (\Z \oplus i \Z)$, hence it defined a holomorphic map
from $X_z$ to the torus $\C/ A^{-1} (\Z \oplus i \Z)$. 

Conversely, if the unfolding $(X_z, \omega_z)$ is a branched cover of the $2$-torus, then $(X_z, \omega_z)$ does not have
the weak mixing property, since non-trivial eigenfunctions for linear toral flows lift to non-trivial eigenfunctions for directional flows
on $(X_z, \omega_z)$.
 \end{proof} 

\begin{cor} \label{cor:ai}  If $k =3, 4, 6$, a polygon $P$ does not have the weak mixing property if and only if it is almost integrable.
\end{cor} 
\begin{proof}  By definition a polygon $P$ is almost integrable if the group $G_P$ generated by reflections with respect to
its edges is an infinite discrete subgroup of the group $E(\R^2)$ of rigid motions. There are only $4$ such groups generated by reflections on the
side of the integrable polygons: rectangles, equilateral triangle, triangle $(\pi/2, \pi/4, \pi/4)$ and triangle $(\pi/2, \pi/3, \pi/6)$.
Let $T_\Delta$ denote the torus corresponding to the integrable polygon $\Delta$. Then $P$ is almost integrable if and only if there exists an integrable polygon $\Delta$ such that the unfolding of $P$ is a branched cover of $T_\Delta$. 
Indeed, is $P$ is almost integrable, then it is drawn on the
graph or ``grid'' generated by reflections of the integrable polygon $\Delta$, hence the unfolding of $P$ is a branched cover of  $T_\Delta$ (see \cite{Gu86}); conversely, if the unfolding of $P$ is a branched cover of a torus $T_\Delta$, it implies that the surface generated by reflections of $P$ is a finite cover of the plane, hence the subgroup generated by such reflections is a discrete subgroup of the group $E(\R^2)$ of rigid motions.  It then follows from Corollary \ref{cor:bct} that if $P$ does not have the weak mixing property, then $P$ is almost integrable.  Conversely, if $P$ is almost integrable then its unfolding is a branched cover of a torus, hence it does not have the weak mixing property (see for instance \cite[\S 3]{Gu84}). 
\end{proof} 

\begin{cor}  If $k =2$, then  $P$ has the weak mixing property if and only if its horizontal side lenghts, as well as its vertical side lengths, are not commensurable.
\end{cor} 
\begin{proof} In this case $T =- Id$ on the tautological plane. If $P$ does not have the weak mixing property, by Theorem \ref{thm:main}  there exists an element $a \re (\omega_z) + b \im (\omega_z) \in H^1(X_z, \Z)\setminus \{0\}$.  Since the edges of the polygon are horizontal or vertical the length of the horizontal/vertical edges are linear combinations with integer coefficients of the half the length of waist curves of horizontal/vertical cylinders.  Since the class $a \re (\omega_z) + b \im (\omega_z)$ is integral $a\not=0$ or $b\not =0$. It follows that  the lengths of horizontal cylinders belong to $a^{-1} \Z$ or the lengths of vertical cylinders belong to  $b^{-1} \Z$, hence  the horizontal or the vertical edges of the polygon are commensurable.

Conversely, if either the  horizontal, or vertical, edges of the polygon have commensurable lengths, it follows that heights of vertical, or horizontal,  cylinders are commensurable, that is, they are integers up to a dilation factor. It is then possible to define a translation map from $X_z$ to the translation circle $\R/ \Z$ by mapping every cylinder to a height segment modulo $\Z$ so that the image of each waist  circle is a single point. Such a map extends continuously to the whole surface $X_z$ since it is the union of horizontal, or vertical, cylinders and all singular leaves at the boundary of cylinders are mapped to the origin of $\R/ \Z$.
By Theorem \ref{thm:main1}, it follows that the polygon $P$ does not have the weak mixing property.
\end{proof}

\begin{ex}  Among triangles, only $(\pi/2, \pi/4, \pi/4)$, $(\pi/2, \pi/3, \pi/6)$ (which are integrable)  and $(2\pi/3, \pi/6, \pi/6)$ 
(which is almost integrable) do not have the weak mixing property. 
\end{ex} 
\begin{proof} By Corollary \ref{cor:gen_case}, a triangle has the weak mixing property unless the least common denominators
of the ration of its angles to $\pi$ is $k=2, 3, 4, 6$.  The case $k=2$ cannot occur for triangles. In the case $k=3$, the triangle
is necessarily equilateral, that is $(\pi/3, \pi/3, \pi/3)$, hence completely integrable. In the case $k=4$, the triangle is necessarily
$(\pi/2, \pi/4, \pi/4)$ which is also integrable. In the case $k=6$, the triangle could be $(\pi/2, \pi/3, \pi/6)$, which is integrable, or
$(2\pi/3, \pi/6, \pi/6)$.  The latter triangle is almost integrable since it is drawn on the lattice generated by the completely integrable
polygon $(\pi/2, \pi/3, \pi/6)$, hence by Corollary \ref{cor:ai} it does not have the weak mixing property.
\end{proof}

\section{Open problems}

We conclude the paper with a list of open problems. An asterisk *  indicates problems
which we believe are beyond the reach of current
methods.

\subsection{Billiards in rational polygons and translation surfaces}
\begin{enumerate}
\item * Prove an explicit upper bound on the Hausdorff dimension of the set of directions that are not weakly mixing whenever
this set has a non-empty complement.  Note, that it is unclear to us how to have an explicit upper bound for the set of directions with an eigenvalue coming from the strong stable direction. We feel that questions like this about the behavior of the cocycle are of interest independent to applications to the weak mixing 
property.
\item Prove polynomial lower bounds for the speed of convergence of Cesaro averages
of correlations  for the vertical flow  (Masur-Veech) typical surface and for almost all directions in Veech surfaces. That is, show that for almost every $(X,\omega)$ there exists an exponent $\sigma>0$ and that there exist Lipschitz functions $f,g \in Lip(X,\omega)$ so that, for all sufficiently large $T>0$,
$$\frac 1 T\int_0^T   \left| \int_X f \circ \phi_t  \cdot g\, d\vert\omega\vert -\int_X f d\vert\omega\vert \int_X g d\vert\omega\vert \right| dt>T^{-\sigma}\,.$$
\item Prove polynomial upper bounds  for the speed of convergence of Cesaro averages
of correlations in almost all directions in non-arithmetic Veech surfaces, that is,
that there exists an exponent $\tau>0$ such that for all smooth functions $f$, $g$ on $X$ there is a constant $C_{f,g}$ such that
$$
\frac 1 T\int_0^T   \left| \int_X f \circ \phi_t  \cdot g \, d\vert\omega\vert -\int_X f d\vert\omega\vert \int_X g d\vert\omega\vert \right| dt  \leq C_{f,g} T^{-\tau}\,.
$$
\item Prove a (hopefully genus independent) upper bound on the Hausdorff dimension of the set of directions that are topologically weakly mixing but not weakly mixing. Note that by Theorem~\ref{thm:main} the translation surfaces with no weakly mixing direction have no topologically weakly mixing direction.
\item * (A. Katok) Prove that the translation flow on the unit tangent bundle to every translation surface $(X,\omega)$ is relatively mixing in the following sense: for any functions $f$, $g\in L^2(X\times \T)$ orthogonal to the subspace of invariant functions (that is, by ergodicity in almost all directions, for any functions $f$, $g$ with zero average on almost every invariant surface in $X \times \T$)  
$$
\int_{X\times \T}  (f\circ \phi_t ) g d \vert \omega\vert d\theta \to  0 \,.
$$
The proof when $(X,\omega)= (\T^2, dz)$ is the flat torus is not difficult. In the higher genus case a result
on typical equidistribution of large circles, relevant to this question, was given by the second author and Hubert \cite{CH}.
\item Construct a translation surface $(X, \omega)$ such that its (forward) Teichm\"uller orbit is recurrent to a compact set, its vertical flow is not weakly mixing and it has an eigenvalue $\alpha$ such that, for all
$\vec{n} \in H^1(\Z)$,
$$
d(g_t \cdot \alpha \im(\omega), g_t\cdot \vec{n}) \to \infty\,.
$$
This will be even more impactful if the example is Birkhoff and Oseledets generic point for the Teichm\"uller flow and the Kontsevich--Zorich cocycle respectively.
\item* Are the translation flows on non-arithmetic Veech surfaces rigid in almost all directions? Alternatively, mild mixing (that is, without rigid factors)? 
\end{enumerate}

\subsection{Billiards in non-rational polygons  }

\begin{enumerate}
\item Find  explicit conditions for the billiard flow of a non-rational polygon to be weakly mixing (on the whole unit tangent bundle).
\item* (Boshernitzan) Does there exist a polygon with minimal billiard flow? (the conjectural answer is negative). Note that one can think of this as a counterpoint to the classic question, does every polygonal billiard have a periodic orbit.  
\item* Does there exist a polygon with mixing billiard flow? (the conjectural answer is yes on the basis of numerical simulations on acute triangles).
\item*
Is the billiard flow in a polygon topologically transitive, or ergodic, or weakly mixing, or mixing for almost all choices of its angles?
(the conjectural answer is positive for weak mixing, hence for ergodicity and topological transitivity. It is more uncertain for mixing).
\end{enumerate}

\bibliographystyle{amsalpha}
\bibliography{bibliography}

\end{document}